\documentclass{article}
\usepackage[margin=1in]{geometry}

\usepackage{microtype}
\usepackage{graphicx}
\usepackage{subfigure}
\usepackage{booktabs}

\usepackage{natbib}
\usepackage{bm}  
\usepackage{algorithm, algorithmic}
\usepackage{graphicx}
\usepackage{multirow,wrapfig}
\usepackage{amsmath, amsthm, amsfonts, amssymb}
\usepackage{booktabs}
\usepackage[ruled, algo2e]{algorithm2e}
\usepackage{xcolor}
\usepackage[normalem]{ulem}
\usepackage{makecell}
\allowdisplaybreaks
\date{}
\usepackage[unicode=true,
            bookmarks=false,
            breaklinks=false,
            pdfborder={0 0 1},
            colorlinks=true,
            citecolor=blue,
            filecolor=blue,
            linkcolor=blue,
            urlcolor=blue]{hyperref}
 \usepackage{url}
\usepackage{amsmath,amssymb,amsthm,mathtools}
\usepackage[capitalize,noabbrev]{cleveref}

\theoremstyle{plain}

\usepackage{microtype}
\usepackage{graphicx}
\usepackage{subfigure}
\usepackage{booktabs} 

\usepackage{hyperref}

\usepackage{amsmath}
\usepackage{amssymb}
\usepackage{mathtools}
\usepackage{amsthm}

\usepackage[capitalize,noabbrev]{cleveref}

\theoremstyle{plain}
\newtheorem{theorem}{Theorem}[section]
\newtheorem{proposition}[theorem]{Proposition}
\newtheorem{lemma}[theorem]{Lemma}

\theoremstyle{definition}
\newtheorem{definition}[theorem]{Definition}
\newtheorem{assumption}[theorem]{Assumption}
\theoremstyle{remark}

\usepackage[english]{babel}
\usepackage[utf8]{inputenc}
\usepackage{amsmath, amsthm, amsfonts, amssymb, bm}
\allowdisplaybreaks
\usepackage{parskip}
\usepackage{graphicx}
\usepackage{setspace}
\usepackage{algorithm}
\usepackage{algorithmic}
\usepackage{natbib}
\usepackage{bbm}
\usepackage[ruled, algo2e]{algorithm2e}
\usepackage{tikz}
\usepackage{xparse}

\usetikzlibrary{decorations.pathreplacing}
\usetikzlibrary{arrows.meta}

\hypersetup{
 colorlinks,citecolor=blue,filecolor=blue,linkcolor=blue,urlcolor=blue}

\newcommand{\VI}{V^{\mathrm{Fluid}}}

\newcommand{\prn}[1]{\left({#1}\right)} 
\newcommand{\prnbig}[1]{\big({#1}\big)} 
\newcommand{\brk}[1]{\left[{#1}\right]} 
\newcommand{\ex}[2]{\mathbb{E}_{#1}\left[#2\right]}

\newcommand{\norm}[1]{\left\|{#1}\right\|}
\newcommand{\Pcal}{\mathcal P}

\newcommand{\Ecal}{\mathcal{E}}
\newcommand{\R}{\mathbb{R}}

\newcommand{\Ical}{\mathcal{I}}

\renewcommand{\P}{\mathbb{P}}
\newcommand{\F}{\mathcal{F}}

\newcommand{\real}{\mathbb{R}}

\newcommand{\Rcal}{\mathcal R}
\newcommand{\N}{\mathcal N}

\newcommand{\st}{\text{s.t.}}

\newcommand{\mix}{\text{\textbf{Hybrid}}}

\definecolor{arc}{RGB}{128,128,255}

\newcommand{\off}{\textbf{off}}
\newcommand{\on}{\textbf{on}}

\newcommand{\KL}[2]{\mathsf{KL}\prnbig{{#1}\,\|\,{#2}}}

\newcommand{\dist}[1]{\mathsf{Dist}\prnbig{{#1}}}
\NewDocumentCommand{\regret}{O{} m}{%
  \mathsf{Regret}^{#1}\prnbig{#2}%
}

\usepackage[textsize=tiny]{todonotes}

\begin{document}
\title{Learning to Price with Resource Constraints: From Full Information to Machine-Learned Prices}
\author{Ruicheng Ao\thanks{R. Ao and D. Simchi-Levi are with the Institute of Data, System and Society at Massachusetts Institute of Technology; emails: {\tt\small \{aorc, dslevi\}@mit.edu}.}\\
Massachusetts Institute of Technology \\
\and
Jiashuo Jiang\thanks{J. Jiang is with Department of Industrial Engineering and Decision Analytics, Hong Kong University of Science and Technology, Hong Kong, China; email: {\tt\small jsjiang@ust.hk}.} \\
Hong Kong University of Science and Technology \\
\and David Simchi-Levi\footnotemark[1]\\
Massachusetts Institute of Technology
}
\maketitle

\begin{abstract}
    We study the dynamic pricing problem with knapsack, addressing the challenge of balancing exploration and exploitation under resource constraints. We introduce three algorithms tailored to different informational settings: a Boundary Attracted Re-solve Method for full information, an online learning algorithm for scenarios with no prior information, and an estimate-then-select re-solve algorithm that leverages machine-learned informed prices with known upper bound of estimation errors. The Boundary Attracted Re-solve Method achieves logarithmic regret without requiring the non-degeneracy condition, while the online learning algorithm attains an optimal $O(\sqrt{T})$ regret. Our estimate-then-select approach bridges the gap between these settings, providing improved regret bounds when reliable offline data is available. Numerical experiments validate the effectiveness and robustness of our algorithms across various scenarios. This work advances the understanding of online resource allocation and dynamic pricing, offering practical solutions adaptable to different informational structures.

\end{abstract}

\section{Introduction}

Dynamic pricing is a classical problem in online learning and decision-making. With the growth of e-commerce, there has been an increasing focus on the development of efficient dynamic pricing policies in the recent literature, see for example \cite{amin2014repeated, javanmard2019dynamic, shah2019semi, cohen2020feature, xu2021logarithmic, bu2022context, xu2022towards}.

Dynamic pricing is challenging due to the unknown relationship between the posted price and the corresponding demand. Usually, one can model the dynamic pricing problem as a bandit problem to balance the exploration of the price-demand relationship and the exploitation of setting the optimal price, where the price option is regarded as an arm and the corresponding revenue (demand multiplied price) is regarded as the reward for the arm. However, different from the classical multi-arm-bandit model where there is a finite number of arms, the price option can be infinite and continuous, which brings significant challenge to the learning part.  

Another challenge to dynamic pricing is the existence of resource constraints. In practice, a product usually enjoys a fixed initial inventory to be sold during the horizon and the pricing decision is made to balance the demand such that the total revenue is maximized subject to the resource constraints. Note that the existence of resource constraints will significantly complicate the problem because the optimal decision is no longer to maximize the revenue but to balance the revenue versus resource consumption in an optimal way. Past work on dynamic pricing mainly focuses on dealing with the learning challenge, but ignores the challenge of resource allocation (e.g. \cite{keskin2014dynamic}). In our paper, we address this fundamental issue by developing near-optimal learning policies for dynamic pricing problems with resource constraints.


\subsection{Preliminaries}
In this paper, we consider an online learning problem with resource constraints, under the context of dynamic pricing. There are $n$ products, $m$ resources and a finite horizon of $T$ discrete time periods. The initial capacities of the resources are denoted by $\bm{C}=(C_1,\dots, C_m)\in\mathbb{R}^+$. For each product $j\in[n]$ and each resource $i\in[m]$, we denote by $a_{ij}$ the amount of resource $i$ required by product $j$, and we denote by $A=[a_{ij}]_{i\in[m],j\in[n]}$ the consumption matrix. 


At the beginning of each period $t$, we denote by $\bm{c}^t=[c_1^t,\dots,c_m^t]$ the remaining capacities of the resources. The decision maker (DM) needs to make an online decision to decide the price for each product, denoted by $\bm{p}^t=(p_1^t, \dots,{\bm{p}}^t_m)\in \Pcal=[L,U]^n$ the price vector, where $L$ (resp. $U$) is a pre-given lower (resp. upper) bound over the product price. The resulting demand is a function over the products' prices and we denote by $f:\Pcal \mapsto\R_+^n$ the demand function. Following classical work on dynamic pricing problems (e.g. \cite{besbes2015surprising}), we adopt the linear demand model by making the following assumption.
\begin{assumption}[Linear demand model]\label{assumption}
There exists parameters $\bm{\alpha}=(\alpha_1, \dots, \alpha_n)\in\R^n$ and $B=[\beta_1,\dots,\beta_n]^\top\in \R^{n\times n}$ such that $f(p) = \alpha+Bp$ for any $\bm{p}\in\Pcal$.
\end{assumption}
Moreover, we make the following assumption on the invertibility and negative-definiteness of the demand model.
\begin{assumption}
    \label{asssump:definite}
    The coefficient matrix $B$ is negative definite, i.e. the minimal eigenvalue $\lambda_{\min}(B+B^\top)>0.$
\end{assumption}
Note that it implies that $B$ is invertible (otherwise we assume there exists $\bm{p}$ such that $B\bm{p} = 0,$ then it follows that ${\bm{p}}^\top (B+B^\top)\bm{p}=0,$ a contradiction).

We assume that there are some noises over the demand side. Therefore, the realized demand, denoted by $\bm{d}^t\in\R_+^n$, under the pricing decision $\bm{p}^t$ is given by $\bm{d}^t=f(\bm{p}^t)+\bm{\epsilon}^t$, where $\bm{\epsilon}^t=(\epsilon_1^t, \dots, \epsilon_n^t)\in \R^n$ are random noises drawn independently from zero-mean sub-Gaussian distributions satisfying $\P(|v^\top\epsilon^t|\ge \lambda)\le 2\exp(-\lambda^2/(2\sigma^2)),\forall t\in[T],\norm{v}_2=1,v\in\R^n$, for some variance parameter $\sigma$. The collected revenue is given by $r^t=(\bm{p}^t)^\top\bm{d}^t$ and the resource consumption is given by $A\bm{d}^t$. The remaining capacities of the resources are then updated by $\bm{c}^{t+1} = \bm{c}^t-A\bm{d}^t$. Note that the DM must adhere to the resource capacity constraints and cannot offer products to the customers when there are no enough remaining capacities. The filtered probability space is denoted by $(\Omega,\F,\{\F^t\}_{t=0}^T,\P)$, where $\mathcal{F}^t$ contains all information by the end of period $t$.


The objective of DM is to maximize the expected revenue over the time horizon $T$ with an online policy $\pi$. The online policy $\pi=(\pi^1,\dots,\pi^T)$ consists of mappings $\pi^t(\bm{c}^t)=\bm{p}^t$ that utilizes only historical information up to time $t$, i.e. $\pi^t$ is measurable on $\F^{t-1}$. The online problems can be formulated as:
\begin{equation}\label{prob:online}
    \begin{aligned}
        \max_\pi \qquad & \Rcal^T(\pi):=\ex{}{\sum_{t=1}^T r^t}=\ex{}{\sum_{t=1}^T(\bm{p}^t)^\top 
 \bm{d}^t}\\ 
        \st \qquad & \bm{d}^t = \bm{\alpha} + B\bm{p}^t +\bm{\epsilon}^t,\\ 
        & A\bm{d}^t \le \bm{c}^t,\forall t\in[T].
    \end{aligned}
\end{equation}

Although the optimal online policy of \eqref{prob:online} can be obtained through dynamic programming (DP), the DP can be computationally intractable due to the curse of dimensionality. An alternative approach is to replace the realized demand by its expectation, resulting in the following fluid model: 
\begin{equation}
    \label{prob:fluid}
    \begin{aligned}
        \max_{\bm{p}\in\Pcal} \qquad & r=\bm{p}^\top \bm{d}\\ 
    \st \qquad & \bm{d} = \bm{\alpha}+B\bm{p},\\ 
    & A\bm{d}\le \frac{\bm{c}}{T}.
    \end{aligned}
\end{equation}
Let $r^\star$ and $\bm{p}^\star$ denote the optimal value and the optimal pricing of the fluid problem \eqref{prob:fluid}, respectively. The fluid model ignores the randomness caused by demand noises and is a deterministic quadratic programming problem that can be solved efficiently. 
The following proposition shows that the fluid model provides an upper bound of \eqref{prob:online}.
\begin{proposition}[\citealt{gallego1994optimal}]
    \label{prop:upper_bound}
    For any online policy $\pi$, we have $Tr^\star\ge \Rcal^T(\pi).$
\end{proposition}
The performance of any policy $\pi$ is measured by the \textit{regret}, which is defined as
\begin{equation}
    \label{eq:regret}
    \regret[T]{\pi} = Tr^*-\Rcal^T(\pi).
\end{equation}

In this paper, we investigate algorithms that reach near-optimal 
performance with low regret.

\subsection{Main Results and Contributions}
We now summarize our main results and main contributions. 

\textbf{Boundary Attracted Re-solve Method}. Our first main result is the development of a novel Boundary Attracted Re-solve Method (Algorithm \ref{alg:resolve}) that achieves a logarithmic regret. Note that the dynamic pricing problem is difficulty even under a full information setting where the price-demand relationship is given, due to the existence of the resource constraints. Efficient policies with a $O(\sqrt{T})$ regret have been developed in the previous literature (see e.g. \cite{balseiro2020dual}). When one seek for an improved logarithmic or better regret, a common assumption has to be made is called \textit{non-degeneracy}, which requires the fluid model \eqref{prob:fluid} to enjoy a unique optimal solution and unique optimal basis (e.g. \cite{li2022online, wang2022constant}). However, as noted in \cite{bumpensanti2020re} and \cite{vera2021bayesian}, degeneracy can indeed happen in practice. In order to relax this non-degeneracy assumption, we develop a new algorithm that resolves the fluid model \eqref{prob:fluid} and further adjust the corresponding optimal solution according to a boundary attracted techniques, i.e., we manually adjust the optimal solution when the remaining resources arrive to the boundary where the degeneracy is about to happen. Such a modification will further push the remaining resources away from the boundary thus avoid the degeneracy issue. We show that our algorithm is able to achieve the state-or-art $O(\log T)$ regret, without any further assumption on the problem instance. 

\textbf{Online Learning with No Prior Information}. We further consider integrating the online learning techniques into our resolving algorithm. We consider the setting where the price-demand relationship is unknown but needs to be learned in an online manner. We show that our algorithm is able to achieve a $O(\sqrt{T})$ regret bound. Most importantly, we show that the $O(\sqrt{T})$ regret is the best we can hope for by establishing a matching lower bound. Therefore, our algorithm achieves the optimal performance in the absence of prior knowledge about the demand parameters. Note that in this setting, the boundary attraction techniques not only serve as a technique to relax the non-degeneracy assumption, but also force exploration over the continuous range of the price options, which leads to our final regret bound. 

\textbf{Online Policies with Informed Price}. We investigate an intermediate setting between the full information setting and the zero information setting. Note that in practice, we usually have a bunch of offline data such that we can deploy some machine learning methods to obtain some estimations/predictions over the price-demand relationship. Such estimations/predictions provide some detailed price-demand information, which differentiates from the zero information setting, but the estimations/predictions can be inaccurate, which differentiates from the full information setting. As for the formulation of the estimations/predictions, we consider the \textit{informed price}, i.e., given a price, we have an estimation/prediction of the expected demand under this price. We explore how to boost the performances of our policies with the help of the informed price. We show that knowing the quality of the informed price indeed matters. To be specific, denote by $\epsilon^0$ the estimation error of the informed price. When the value of $\epsilon^0$ is not known, we develop a lower bound showing that one cannot achieve a worst-case regret better than $O(\sqrt{T})$. On the other hand, when the value of $\epsilon^0$ is given, we develop a new estimate-then-select re-solve algorithm (Algorithm \ref{alg:resolve_learn_incumbent}) that leverages the machine-learned informed price-demand pair and we show that our algorithm achieves a regret bound of $O(\min\{\sqrt{T}, (\epsilon^0)^2\cdot T\}+\log T)$, which is further justified to be of optimal order with a matching lower bound. As we can see, as long as the informed price is accurate enough, i.e., $\epsilon^0\leq1/\sqrt{T}$, we recover the $O(\log T)$ regret under the full information setting. Our results bridge the gap between the full information and no information settings, providing improved regret bounds when reliable offline data is available.



    
    
\subsection{Other Related Literature}
Our paper is related to three streamlines of literature: (i) bandit with knapsack; (ii) online resource allocation with fluid approximation; (iii) dynamic pricing with offline data and misspecification.


\noindent\textbf{Bandits with Knapsacks}. The bandits with knapsacks (BwK) framework, introduced by \cite{agrawal2016linear}, is closely related to our dynamic pricing problem. In CwK, rewards are disclosed to the agent following the decision-making process. Two main methodologies have been developed to address this problem. The first method focuses on selecting the optimal probabilistic strategy from a predefined policy set \cite{badanidiyuru2014resourceful}, with \cite{agrawal2016efficient} employing this technique to achieve an \(O(\sqrt{T})\) regret bound and \cite{liu2022non} extends to the non-stationary setting. This approach is rooted in contextual bandit theory \cite{dudik2011efficient,badanidiyuru2014resourceful} and utilizes a cost-sensitive classification oracle to ensure computational efficiency. Alternatively, the second methodology approaches the CBwK problem through the Lagrangian dual framework. It applies a dual update mechanism that transforms the BwK problem into an online convex optimization (OCO) problem. Specifically, studies such as \cite{agrawal2016linear,sankararaman2021bandits,sivakumar2022smoothed,liu2022online} consider a linear model. This approach integrates linear function estimation techniques \cite{abbasi2011improved,auer2002using,sivakumar2020structured,elmachtoub2022smart,kumar2022non,zhang2024piecewisestationary} with OCO methods to achieve sub-linear regret. More recently, \cite{chen2024contextual} considers contextual decision-making with knapsack problem and reaches constant regret beyond worst case.

\textbf{Online resource allocation with fluid approximation.} A line of research related to our model and analysis framework investigates the near-optimal online algorithms comparing with the offline benchmark in online resource allocation problems \citep{reiman2008asymptotically,jasin2012re,jasin2013analysis,ferreira2018online,bumpensanti2020re,banerjee2020uniform,vera2021bayesian,wang2022constant,jiang2022degeneracy,jaillet2024should,ao2024online,ao2024two}. \cite{wang2022constant} proves that fluid approximation based re-solve methods can achieve near-optimal logarithmic regret compared with the fluid approximation benchmark in pricing-based revenue management problem under a non-degeneracy condition. \cite{broder2012dynamic,chen2021joint,Chen2022} consider a maximum-likelihood estimator and adopts an explore-then-exploit policy. \cite{broder2012dynamic,chen2021joint} proves a logarithmic regret under a separate condition.

\textbf{Pricing with offline data and misspecification}
\cite{keskin2014dynamic} considers an ``incumbent price" setting that assumes existence of an accurate data point. \cite{bu2020online} discuss
a phase transition phenomenon to quantify how useful the offline data are for an online
pricing problem. \cite{li2021unifying} and \cite{wang2023measuring} also consider the settings when offline dataset is available when doing online decision. \cite{ferreira2018online} consider the model misspecification of online pricing problem. More recently, \cite{wang2024online} discuss the presence of polluted offline data, which is more relevant to our paper.

\textbf{Notation}. 
For $x\in\real$, denote $\lceil x\rceil$ as the smallest integer not smaller than $x$ and $\lfloor x\rfloor$ as the largest integer not greater than $x$. 
Denote $x_+=\max\{x, 0\}$. For set $S$, denote $|S|$ as its cardinality. 
Additionally, denote $d_{\max} = \max_{\bm{p}\in\P}\norm{f(\bm{p})}_2$ as the norm of maximal feasible demand.

\section{Algorithm and Logarithmic Regret with Full Information}
In this section, we assume that we have full information over the demand function $f$ beforehand and we present our main algorithmic framework. Given the remaining capacities $\bm{c}^t$ at period $t$, the fluid model can be written as
\begin{equation}
    \label{prob:resolve}
    \begin{aligned}
        \VI_t(\bm{c}^t)=\max_{\bm{p}\in\Pcal}\qquad & r = \bm{p}^\top \bm{d}\\ 
    \st\qquad& \bm{d} = \bm{\alpha} + B\bm{p},\\ 
    &A\bm{d}\le \frac{\bm{c}^t}{T-t+1},
    \end{aligned}
\end{equation}
A natural idea is to resolve the fluid model $\VI_t$ in \eqref{prob:resolve} to obtain the pricing decision $\bm{p}^t$. However, as noted in a series of work (e.g. \cite{jasin2014reoptimization} and \cite{wang2022constant}), such resolving methods achieve near-optimal performances only when the non-degeneracy condition is guaranteed, i.e., the optimization problem \eqref{prob:resolve} preserves a unique optimal basis during the horizon. In order to relax this non-degeneracy assumption, we introduce a ``boundary attraction'' technique that absorbs the price to the upper bound $U$ when the corresponding expected demand is sufficiently small. Our formal algorithm is presented in Algorithm \ref{alg:resolve}. We now show that our Algorithm \ref{alg:resolve} achieves a logarithmic regret with no further condition.

  \begin{algorithm}[tb]
	\caption{Boundary attracted re-solve method: known parameters}
	\begin{algorithmic}[1]\label{alg:resolve}
		\STATE{\textbf{Input:} Initial inventory level $\bm{c}^1=\bm{C}$, consumption matrix $A$, demand function $f(\bm{p})=\bm{\alpha}+B\bm{p}$ and rounding parameter $\zeta$.}
		\FOR{$t=1,2,\dots,T$} 
                \STATE{Obtain $(\bm{p}^{\pi,t},\bm{d}^{\pi,t})$ as one optimal solution of $\VI_t(\bm{c}^t)$ in \eqref{prob:resolve}.}
                \STATE{Set targeted demand $\tilde{\bm{d}}^t$ according to
                \begin{equation*}
                    \tilde{d}_i^t = \left\{\begin{aligned}
                        d_i^{\pi,t}, &\quad d_i^{\pi,t}\ge \zeta(T-t+1)^{-1/2},\\ 
                        0, &\quad d_i^{\pi,t}< \zeta(T-t+1)^{-1/2}.
                    \end{aligned}\right.
                \end{equation*}}
			\STATE{Set actual price $\bm{p}^t$ such that $f(\bm{p}^t)=\tilde{\bm{d}}^t$.}
			\STATE{Observe realized demand $\hat{\bm{d}}^t=f(\bm{p}^t)+\bm{\epsilon}^t$.}
			\STATE{Update inventory level $\bm{c}^{t+1}=\bm{c}^t-A\hat{\bm{d}}^t$.
			}
		\ENDFOR
	\end{algorithmic}
\end{algorithm}

\begin{theorem}
    \label{thm:full_information}
    Let $\pi$ be given by Algorithm \ref{alg:resolve}, assume that $\zeta\ge 4\sigma^2,$ then we have $\regret[T]{\pi} =  O\prn{\zeta^2n^2\norm{B^{-1}}_2\log T}.$
\end{theorem}
We sketch our proof here, while the details are left to Appendix \ref{appendix:resolve}. Our steps can be splitted into three steps: \textbf{Step I}: we introduce the so-called ``hybrid" policy $\mix^t$ that adopts our Algorithm \ref{alg:resolve} up to time $t$ and get the remain revenue by solving \eqref{prob:resolve} directly without noise. We can then decouple the regret into single-step difference of the form 
\begin{align*}
&\regret[T]{\pi} \\&= \ex{}{\sum_{t=1}^T\Rcal^T(\mix^t,\F^T)-\Rcal^T(\mix^{t+1},\F^T)},
\end{align*}
where $\Rcal^T(\mix^t,\F^T)$ denotes the total revenue of the hybrid policy defined above under realized sample paht $\F^T.$ By the decomposition above, we can give upper bound to the total revenue by analyzing each single-step difference and sum them up. \textbf{Step II}: Upon observing noise $\bm{\epsilon}^t$ at time $t$, in order to control the single-step difference, we apply a shift of $\bm{\epsilon}^t/(T-t)$ to the demands in the remain $T-t$ time horizons. Specifically, since $\mix^{t+1}$ will take optimal fluid optimal for the remain periods $t+1,\dots,T$ without noise, it would not generate higher revenue by substituting it with any feasible pricing-demand pair for the remain $T-t$ periods. Particularly, if $d_i^{\pi,t}>\epsilon_i^t/(T-t),$ for all $i\in[n],$ (which will hold with high probability if all $d_i^{\pi,t}$ are over the attraction threshold) then we can substitute the fluid solution in $\mix^{t+1}$ with $\bm{d}^{\pi,t}-\bm{\epsilon}^t/(T-t),$ which is feasible since $(T-t)A(\bm{d}^{\pi,t}-\bm{\epsilon}^t/(T-t)) + (\bm{d}^{\pi,t}+\bm{\epsilon}^t) \le \bm{c}^t.$ This is the key step to hedge the first-order error caused by the noise with accumulating errors of order $1/(T-t).$ \textbf{Step III}: When $d_i^{\pi,t}$ is attracted to the boundary, we can simply replace the previous $d_i^{\pi,t}-\epsilon_i^t/(T-t)$ by $(T-t+1)d_i^{\pi,t}/(T-t),$ which is still a feasible solution to the fluid problem and it will generate an error of order $1/(T-t+1)$. By combining the steps above, we can conclude that, in each step, the error grows in order of $1/(T-t+1),$ which will lead to a total regret of $O(\log T).$

Note that we donot require the non-degeneracy condition proposed in previous literatures \cite{jasin2014reoptimization,wang2022constant}.  Compared with the results in \cite{jasin2014reoptimization} and \cite{wang2022constant}, the novelty of our Algorithm \ref{alg:resolve} can be summarized in three points: (i) We donot apply the first-order correction step that takes $f({\bm{p}}^t) = f({\bm{p}}^{t-1})-\frac{\bm{\epsilon}^t}{T-t+1}$, but use the accurate optimal solution to \eqref{prob:resolve} instead; (ii) We apply a rounding threshold $\zeta(T-t+1)^{-1/2}$ in order to hedge the risk of ``over-estimation" of the demand in the presence of noises; (iii) The previous work needs non-degeneracy of optimal solution to \eqref{prob:fluid} to guarantee that the first-order correction will generate an optimal solution of the noised version. With a novel single-step difference technique that was first applied in quantity-based revenue management \cite{vera2021bayesian,jiang2022degeneracy}, we can remove the non-degeneracy assumption in previous work.

\section{Algorithm and Regret with No Information}
In this section, we consider a no information setting where the demand function $f(\cdot)$ is unknown and need to be learned in an online fashion as we obtain more and more observations over the realized demand. Under Assumption \ref{assumption}, the learning task can be summarized as learning the parameters $\bm{\alpha}$ and $B$. The most natural way is to adopt linear regression to learn the parameters. To be more specific, at each period $t$, given the data points $(\bm{p}^1,\bm{d}^1),\dots,(\bm{p}^{t-1},\bm{d}^{t-1})$, we define 
\begin{equation}\label{eqn:DefDandP}
\begin{aligned}
    D_j^t &:= \sum_{s=1}^{t-1}[d_j^s; d_j^s\cdot\bm{p}^s]^\top,~\forall j\in[n], \\
    P^t &:=  \begin{bmatrix}
        t-1 & \sum_{s=1}^{t-1}(\bm{p}^s)^\top\\
        \sum_{s=1}^{t-1}\bm{p}^s & \sum_{s=1}^{t-1}\bm{p}^s(\bm{p}^s)^\top
    \end{bmatrix}.
\end{aligned}
\end{equation}
Then, the linear regression can be done by computing
\begin{equation}
\label{eq:linear_regression}
\begin{aligned}
         \begin{bmatrix}
     \hat\alpha_j^t\\ 
     \hat{\bm{\beta}}_j^t
 \end{bmatrix} 
&= (P^t)^{\dag} D_j^t \\ &= \begin{bmatrix}
     \alpha_j\\ 
     \bm{\beta}_j
 \end{bmatrix} + (P^t)^{\dag}\begin{bmatrix}
     \sum_{s=1}^{t-1}\epsilon_j^s\\
     \sum_{s=1}^{t-1}\epsilon_j^s\cdot\bm{p}^s
 \end{bmatrix},\quad\forall j\in[n],
\end{aligned}
\end{equation}
where $(P^t)^{\dag}$ represents the pseudo-inverse of the matrix $P^t$. As we can see from \eqref{eq:linear_regression}, the estimation error (the gap between $(\hat{\bm{\alpha}}^t, \hat{B}^t)$ and $(\bm{\alpha}, B)$) scales linearly in the smallest eigenvalue of the matrix $P^t$, which can be lower bounded by the variance of the historical prices $\bm{p}^1,\dots,\bm{p}^{t-1}$. Indeed, motivated by \cite{keskin2014dynamic}, denoting $\overline{\bm{p}}^t=t^{-1}\sum_{s=1}^t\bm{p}^s,$ one can show that 
\begin{align}\label{eq:fisher}
    &\lambda_{\min}(P^t)
    \\&\ge \frac{1}{n(1+2U^2)}\sum_{s=1}^{n\lfloor t/n\rfloor}(1-\frac{1}{s})\big\|\bm{p}^s-\overline{\bm{p}}^{s-1}-{\bm{X}}^{\lceil s/n\rceil}\big\|_2^2
\end{align}
for some fixed anchor $\bm{X}^{k},k=0,1,\dots,$ such that $\bm{p}^s-\overline{\bm{p}}^{s-1}-{\bm{X}}^{k},s=kn+1,\dots,kn+n$ form an orthogonal basis of $\R^n$. 
Therefore, in order to reduce the estimation error and guarantee sufficient exploration, instead of directly using the optimal solution of the fluid model $\VI$, we add some perturbation on the pricing strategy $\bm{p}^t$ to get a higher variance. Specifically, at time $t$, we can do this by setting the desired price as
\begin{align*}
    \tilde{\bm{p}}^t = \overline{\bm{p}}^{t-1} + (\tilde{\bm{p}}^{k}-\overline{\bm{p}}^{kn}) + \sigma_0t^{-1/4}e_{t-kn}
\end{align*}
for $k = \lfloor (t-1)/n\rfloor$ and $\tilde p^k$ is the solution to \eqref{prob:resolve} by replacing the parameters with estimated ones $\hat{\bm{\alpha}}^{kn+1},\hat B^{kn+1}$ and inventory level $\bm{c}^{kn+1}.$ As a result, the term $(\tilde{\bm{p}}^{k}-\overline{\bm{p}}^{kn})$ performs as a momemtum forward to the re-solve solution and an anchor for time steps $kn+1,\dots,kn+n,$ hence balancing the exploitation and the exploration.

Our algorithm is still developed based on the Boundary attracted resovling algorithm in Algorithm \ref{alg:resolve}. However, different from the full information case, we make modifications to incorporate the learning on demand parameters: (i) we apply linear regression on historical data to fit the model parameters and substitute them into the re-solve problem \eqref{prob:resolve} at each time; (ii) we make a perturbation on the optimal solution to \eqref{prob:resolve} to enforce exploration and make sure the data variance will scale up moderately. 
  \begin{algorithm}[tb]
	\caption{Periodic-review re-solve with parameter learning}
	\begin{algorithmic}[1]\label{alg:resolve_learn}
		\STATE{\textbf{Input:} Initial capacity $\bm{c}^1=\bm{C},$ consumption matrix $A$, variance scale $\sigma_0$, rounding threshold $\zeta$.}
            \FOR{$t=1,2,\dots,n$}
                \STATE{Uniformly sample $\bm{p}^t$ from $\Pcal$.}
            \ENDFOR
            \STATE{Calculate average price $\overline{\bm{p}}^n=\frac{1}{n}\sum_{t=1}^n\bm{p}^t.$}
		\FOR{$t=n+1,n+2,\dots,T$} 
                \IF{$\mod(t,n)=1$}
                    \STATE{Set $k=(t-1)/n$.}
                    \STATE{Calculate coefficients $\hat{\bm{\alpha}}^{kn+1},\hat{B}^{kn+1}$ according to \eqref{eq:linear_regression}.}
                    \STATE{Calculate $\tilde{\bm{p}}^{k}$ by solving \eqref{prob:resolve} with parameters $\hat{\bm{\alpha}}^{kn+1},\hat{B}^{kn+1}$ and remaining capacities $\bm{c}^{kn+1}$.}
                \ENDIF
			\STATE{Calculate price ${\bm{p}}^t$ as 
            \begin{equation*}
{\bm{p}}^t = \overline{\bm{p}}^{t-1} + (\tilde{\bm{p}}^{k}-\overline{\bm{p}}^{kn}) + \sigma_0t^{-1/4}e_{t-kn}
            \end{equation*}}
            \vspace{-0.1in}
                \STATE{Calculate the corresponding estimated demand $\tilde{\bm{d}}^t= \hat{\bm{\alpha}}^{kn+1}+\hat B^{kn+1} {\bm{p}}^t$.}
                \STATE{Set indices set as $\Ical_r^t=\{i\in[n]:\tilde d_i^t\le \zeta((T-t+1)^{-1/4}+t^{-1/4})\}$. }
                \STATE{Observe actual demand $\hat d^t$ corresponding to price $p^t.$}
                \STATE{Reject demands of type $i\in\Ical_r^t$ and update accepted demands as $\hat d^t$ accordingly.}
			\STATE{Update inventory level ${\bm{c}}^{t+1}={\bm{c}}^t-A\hat{\bm{d}}^t$ with accepted demands $\hat{\bm{d}}^t$.
			}
                \STATE{Update average price $\overline{\bm{p}}^t = \frac{(t-1)}{t}\overline{p}^{t-1}+\frac{1}{t}{\bm{p}}^t.$}
		\ENDFOR
	\end{algorithmic}
\end{algorithm}
Our algorithm is formally presented in Algorithm \ref{alg:resolve_learn} and the regret is bounded below.
\begin{theorem}
    \label{thm:regret_no_information}
        Let $\pi$ be given by Algorithm \ref{alg:resolve_learn}, assume that $\zeta\ge Cn^{5/4}\log^{3/2}n\sigma_0\sqrt{\sigma}\log T$ for some constant $C$, we have $\regret[T]{\pi} = O\prn{(\zeta^2+\norm{B^{-1}}_2)\sqrt{T}}.$ 
\end{theorem}
We sketch the proof for Theorem \ref{thm:regret_no_information} here, while the details are left in Appendix \ref{appendix:learn}. WLOG we assume $T=nT'$ for some integer $T'>0$. The first step is similar to the proof of Theorem \ref{thm:full_information}: we split the regret as
\begin{align*}
    &\regret[T]{\pi}\\ 
    &= \ex{}{\sum_{k=1}^{T'}\Rcal^T(\mix^k,\F^T)-\Rcal^T(\mix^{k+1},\F^T)},
\end{align*}
where $\R^T(\mix^k,\F^T)$ denotes the total revenue of the hybrid policy defined as using Algorithm \ref{alg:resolve_learn} up to time $kn$ and getting the remain revenue by solving \eqref{prob:resolve} directly without noise. Now different from the proof of Theorem \ref{thm:full_information}, we donot have the accurate parameters as well as the fluid optimal solutions. As a result, we need to give bound to the estimation error of $\hat{\bm{\alpha}}^{kn+1},\hat B^{kn+1}$ with respect to the true parameters $\bm{\alpha},B,$ as well as the corresponding estimated solutions $p^t,f(p^t)$. To achieve this target, we briefly introduce the continuity property in strongly convex constrained optimization problem in objective function. 
\begin{lemma}[Prop 4.32, \cite{bonnans2013perturbation}]
    Suppose the constrained optimization problem \eqref{prob:resolve} satisfies the second-order growth condition $\bm{p}^\top\bm{d}-(\bm{p}^{\pi,t})^\top\bm{d}^{\pi,t}\le \kappa\dist{\bm{d},D^{\pi,t}}$ for some constant $\kappa>0$ and feasible solution $\bm{p}^{\pi,t})^\top\bm{d}^{\pi,t}$ to \eqref{prob:resolve}, where $D^{\pi,t}=\{\bm{d}^{\pi,t}:\bm{d}^{\pi,t},\bm{p}^{\pi,t} \text{ solve \eqref{prob:resolve}}\}$. Then we have $\dist{\hat{\bm{d}},D^{\pi,t}}\le C\kappa^{-1}(\norm{B-\hat B}_2+\norm{\bm{\alpha}-\hat{\bm{\alpha}}}_2)$ for some constant $C>0$, where $\hat B,\hat{\bm{\alpha}}$ lies in some neighbor of $B,\bm{\alpha}$ depending on $\lambda_{\min}(B+B^\top)$ and $\hat{\bm{d}}$ is the solution to \eqref{prob:resolve} by replacing $\bm{\alpha},B$ with $\hat{\bm{\alpha}},\hat B.$
\end{lemma}
Now by applying the above lemma, we can bound the estimation error of solution $p^t,f(p^t)$ based on estimation error of the parameters $\hat{\bm{\alpha}}^{kn+1},\hat B^{kn+1}$, which is achievable by \eqref{eq:fisher}. By combining the estimation error and the error caused by noises as in the proof of Theorem \ref{thm:full_information}, we can get the final bound in Theorem \ref{thm:regret_no_information}.

Even in the unconstrained case, it is well-known that a $\sqrt{T}$ lower bound of worst-case regret is inevitable \citep{keskin2014dynamic, Chen2022}. Formally, we have
\begin{lemma}[\citealt{keskin2014dynamic}]\label{lem:lower_bound}
    There exists a finite positive constant $c>0$ such that $\regret[T]{\pi}\ge c\sqrt{T}$ for any online policy $\pi$. 
\end{lemma}
Therefore, we know that the regret bound presented in Theorem \ref{thm:regret_no_information} is of optimal order.

\section{Improvement with Machine-learned Informed Price}
In previous sections, we investigate the full information setting and the zero information setting. In practice, however, we usually encounter an intermediate setting where there can be some offline dataset from which we can form an estimation of the demand function parameters using some machine learning algorithms. Therefore, in this section, we explore how some pre-given machine-learned estimations can be helpful in improving our algorithm performance, which also builds a bridge between the full information setting and the zero information setting. 

We consider the following formulation of the machine-learned estimations, which we call \textit{informed price}. In real-life scenarios, the sellers often set a fixed price to collect demand information \citep{grossman1976information}. Such price is called informed price since it is most reflective of all shared present or past market information. Formally, we assume that we have a pre-given price-demand pair $(\bm{p}^0,\bm{d}^0)$. In practice, the seller can first set the pricing decision to be $\bm{p}^0$ and collect the corresponding demand information to form an estimation over the expected demand, where the estimation is denoted by $\bm{d}^0$. Note that due to some estimation error, $\bm{d}^0$ can be different from $f(\bm{p}^0)$, which represents a misspecification of model parameters. 


We explore how to use $(\bm{p}^0, \bm{d}^0)$ to better make the online pricing decision. We denote by $\epsilon^0$ an upper bound of the estimation error, i.e., we define $\epsilon^0$ such that $\|\bm{d}^0 - (\bm{\alpha}+B\bm{p}^0)\|_2\leq\epsilon^0$. 
We argue that it is indeed important to know the quality of the informed price-demand pair $(\bm{p}^0, \bm{d}^0)$, i.e., an upper bound $\epsilon^0$ on the estimation error. Instead, if we do not have the upper bound $\epsilon^0$, we show in the next proposition that one cannot do better than the no information setting even when the informed price-demand pair $(\bm{p}^0, \bm{d}^0)$ is given (but the estimation error is ungiven).



\begin{proposition}
    \label{prop:impossible}
    There exists two groups of parameters $(\alpha,B)$ and $(\alpha',B')$ such that for any admissible policy $\pi$ without knowledge of the accuracy, if $\pi$ has regret $O(T^\gamma)$ on model parameterized by $(\alpha,B)$ for some $\gamma\in(0,1)$, then it has regret $\Omega(T^{1-\gamma})$ on model $(\alpha',B')$.
\end{proposition}
Proposition \ref{prop:impossible} shows that there is no online policy that can do uniformly better than the zero information setting when the informed price-demand pair is given but an upper bound of the estimation error is ungiven. Note that Proposition \ref{prop:impossible} implies a $\Omega(\sqrt{T})$ regret lower bound for any online policy when we have no knowledge of $\epsilon^0$. On the other hand, when $\epsilon^0$ is not given, we know that the best we can do is to treat the problem under the zero information setting, even when $(\bm{p}^0, \bm{d}^0)$ is given, and apply Algorithm \ref{alg:resolve_learn}.


We now consider how to utilize the informed price-demand pair $(\bm{p}^0, \bm{d}^0)$ to enhance the algorithm performance, assuming that we are given an upper bound $\epsilon^0$ on the estimation error.
Note that by utilizing the informed price-demand pair $(\bm{p}^0,\bm{d}^0)$, the estimation of the linear coefficient $B$ can be fitted by
\begin{equation}
    \label{eq:regression_informed}
    \hat{B}^t = \brk{\sum_{s=1}^{t-1}(\bm{p}^s-\bm{p}^0)(\bm{p}^s-\bm{p}^0)^\top}^{\dag}\sum_{s=1}^{t-1}(\bm{d}^s-\bm{d}^0)(\bm{p}^s-\bm{p}^0)^\top.
\end{equation}
Then the estimated re-solve problem can be formulated as:
\begin{equation}
    \label{prob:resolve_incumbent}
    \begin{aligned}
     \hat{V}_t^{\mathrm{Fluid}}(\bm{c}^t)=   \max_{\bm{p}\in\Pcal}\qquad & r = \bm{p}^\top \bm{d}\\ 
    \st\qquad& \bm{d} = \bm{d}^0 + B(\bm{p}-\bm{p}^0),\\ 
    &A\bm{d}\le \frac{\bm{c}^t}{T-t+1},
    \end{aligned}
\end{equation}
We make the following assumption on the informed price-demand pair.
\begin{assumption}
    \label{asssump:informed}
    We have $\bm{d}^0 > \bm{d}^\star,$ where $\bm{p}^\star,\bm{d}^\star$ are solutions to \eqref{prob:fluid} with infinite resources $\bm{c} = +\infty.$
\end{assumption}
In other word, the informed price-demand pair can be considered as being obtained through discounts or promotions. Note that in this case, different from the lower bound in \eqref{eq:fisher}, the informed price-demand pair can be considered as replications of infinite data points, hence enlarging the variance of the dataset for estimation. As a result, we can expect that setting lower perturbation can yield higher estimation accuracy than the no information case.

We propose an estimate-then-select re-solve algorithm, which decides whether to trust the prediction oracle or not at the beginning, based on the value of $\epsilon^0$. Our algorithm is formally presented in Algorithm \ref{alg:resolve_learn_incumbent}.
  \begin{algorithm}[tb]
	\caption{Estimate-then-select re-solve with parameter learning}
	\begin{algorithmic}[1]\label{alg:resolve_learn_incumbent}
		\STATE{\textbf{Input:} Inventory level ${\bm{c}}^1=\bm{c},$ cost matrix $A$, upper bound $\epsilon^0$ on misspecification, tolerance $\rho$, variance scale $\sigma_0$, rounding threshold $\zeta$.}
            \IF{$(\epsilon^0)^2T > \rho\sqrt{T}$}
                \STATE{Change to Algortihm \ref{alg:resolve_learn}.}
            \ENDIF
		\FOR{$t=1,2,\dots,T$} 
                \STATE{Set $l = \mod(t,n)$.}
                \STATE{Calculate coefficients $\hat B^t$ according to \eqref{eq:regression_informed}.}
                \STATE{Calculate $\tilde{\bm{p}}^t$ by solving \eqref{prob:resolve_incumbent} with estimated parameters $\hat B^t$, inventory level ${\bm{c}}^t$.}
                \STATE{
                    Update $\tilde{\bm{p}}^t$ as 
                    \begin{equation*}
                        \tilde{\bm{p}}^t \leftarrow \tilde{\bm{p}}^t + \sigma_0 \text{sgn}(\tilde{\bm{p}}^t-{\bm{p}}^0)t^{-1/4}e_{l}.
                    \end{equation*}
                }
                \STATE{Calculate the corresponding $\tilde{\bm{d}}^t$ with estimated parameters $\tilde{\bm{d}}^t = \bm{d}^0+\hat B^t(\tilde{\bm{p}}^t-\bm{p}^0)$.}
                \STATE{Set the indices set as $\Ical_r^t=\{i\in[n]:\tilde d_i^t\le \zeta((T-t+1)^{-1/2}+t^{-1/2})\}$}
			\STATE{Observe demand $\hat{\bm{d}}^t=f(\tilde{\bm{p}}^t)+\epsilon^t$.}
                \STATE{Reject demands of types $i\in\Ical_r^t$ and update accepted demands $\hat{\bm{d}}^t$}
			\STATE{Update inventory level ${\bm{c}}^{t+1}={\bm{c}}^t-A\hat{\bm{d}}^t.$
			}
		\ENDFOR
	\end{algorithmic}
\end{algorithm}
The regret guarantee of Algorithm \ref{alg:resolve_learn_incumbent} can be bounded in the following theorem.
\begin{theorem}
    \label{thm:incumbent}
    Let $\pi$ be given by Algorithm \ref{alg:resolve_learn_incumbent}, we have $\regret[T]{\pi} = O(\min\{\rho\sqrt{T},(\epsilon^0)^2T+C'\log T\}),$ where $C' = \sigma_0d_{\max}\norm{B^{-1}}_2+n^2\sigma^2\norm{B^{-1}}_2+\zeta^2.$
\end{theorem}

The proof is similar to the proof of Theorem \ref{thm:regret_no_information} except that we give bound to the estimation error in presence of the informed price-demand pair $\bm{p}^0,\bm{d}^0$. Details are left in Appendix \ref{appendix:incumbent}.

When $\epsilon^0$ is large such that $\epsilon^0\geq\Omega(1/\sqrt{T})$, the regret bound in Theorem \ref{thm:incumbent} recovers the regret under the zero information setting presented in Theorem \ref{thm:regret_no_information}. Instead, when the estimation is accurate such that $\epsilon^0=0$, the regret bound in Theorem \ref{thm:incumbent} recovers the regret under the full information setting presented in Theorem \ref{thm:full_information}. On the other hand, we can also show in the following theorem that the regret in Theorem \ref{thm:incumbent} is of optimal order when an upper bound $\epsilon^0$ is given, which shows the tightness of our results.

\begin{proposition}
    \label{prop:impossible_incumbent}
    There exists two instances $(\alpha,B)$ and $(\alpha',B')$ such that for any admissible policy $\pi$ without knowledge of the accuracy $\epsilon^0$ of the informed demand-price pair, if $\pi$ has regret of $O(\min\{\rho\sqrt{T},(\epsilon^0)^2T\})$ on instance $(\alpha,B)$, then it has regret $\Omega(\max\{\rho\sqrt{T},(\epsilon^0)^2T\}).$
\end{proposition}

\section{Numerical Validations}
In this section, we use numerical experiments to verify our analysis. We perform our simulation experiments in the following three settings: (i) full information; (ii) no information; (iii) informed price. We adopt $A = [1,1],\alpha = [8,6]^\top, B = [-1/2, -1/5;-1/5,-1/2]$ and degenerate budget $c = A{\bm{d}}^{\star},$ with ${\bm{d}}^{\star}$ is obtained by solving $\max_d{\bm{d}}^\top B^{-1}(d-\alpha).$ We set the noises as i.i.d. Gaussian noises $\mathcal N(0,1),$ perturbation scale $\sigma_0=1,$ rounding threshold $\zeta=1,$ misspecification error $\epsilon^0 = T^{-1/2}.$

\begin{figure}[h!]
\vskip 0.2in
    \begin{center}
\centerline{\includegraphics[width=0.8\linewidth]{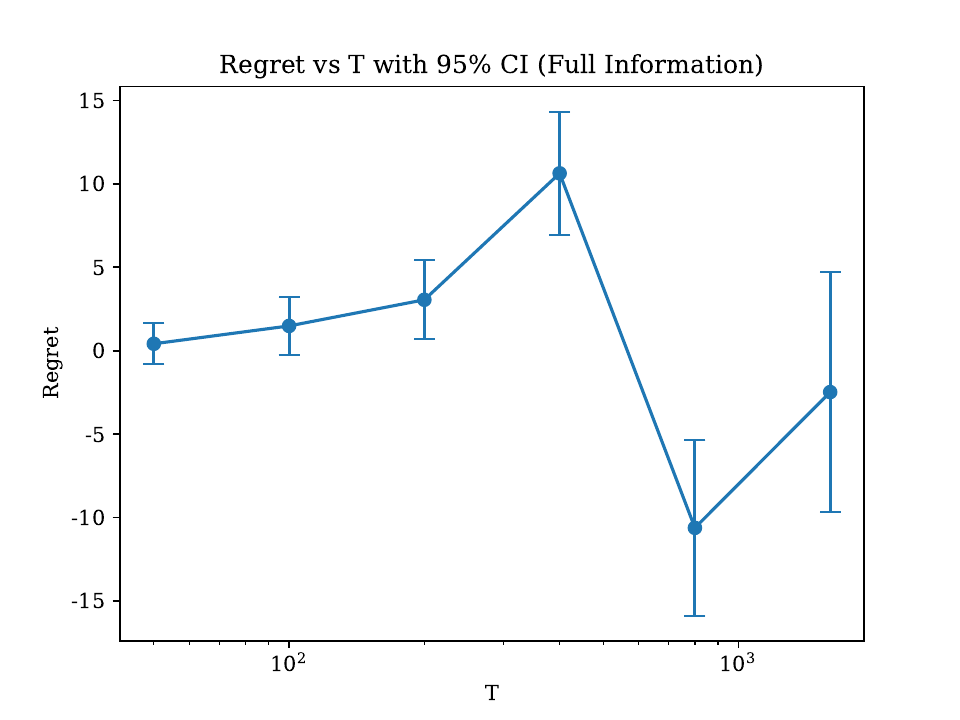}}
\caption{Regret of our algorithms under different number of time horizon $T$ with full information\label{fig:comparison1}.}
\end{center}
\vskip -0.2in
\end{figure}

\begin{figure}[h!]
\vskip 0.2in
    \begin{center}
\centerline{\includegraphics[width=0.8\linewidth]{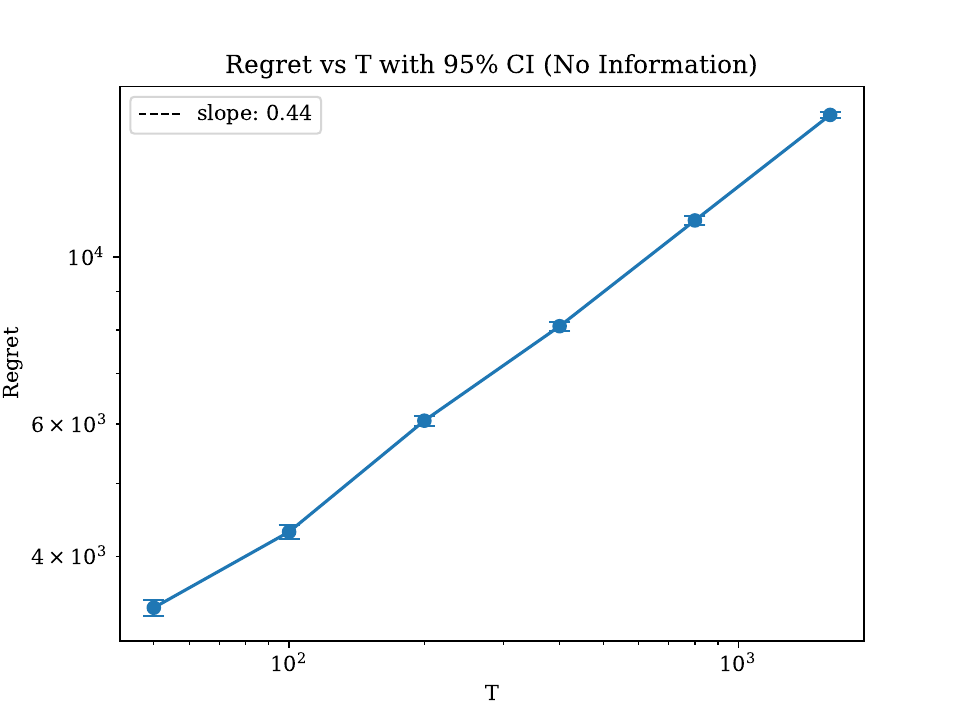}}
\caption{Regret of our algorithms under different number of time horizon $T$ with no information ($\log T$-$\log\text{Regret}$ 
 plot)\label{fig:comparison2}.}
\end{center}
\vskip -0.2in
\end{figure}

\begin{figure}[ht]
\vskip 0.2in
    \begin{center}
\centerline{\includegraphics[width=0.8\linewidth]{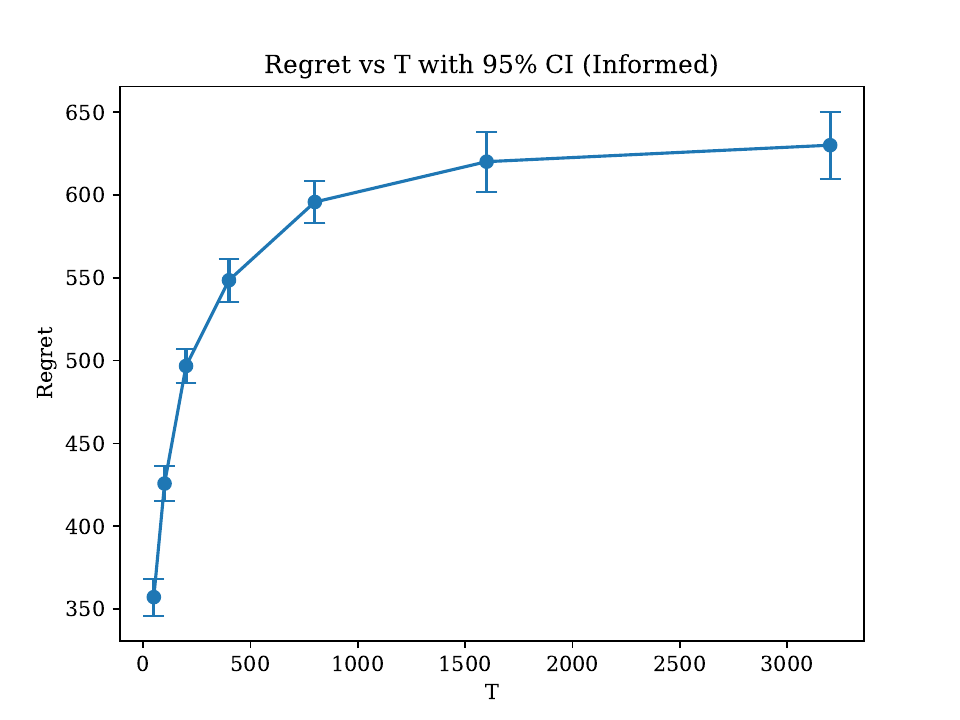}}
\caption{Regret of our algorithms under different number of time horizon $T$ with informed price\label{fig:comparison3}.}
\end{center}
\vskip -0.2in
\end{figure}

Figure \ref{fig:comparison1}, \ref{fig:comparison2} and \ref{fig:comparison3} illustrate the relationship between the regret and the number of time horizon $T$ under the three settings. We set the horizon $T$ to be $\{50,100,200,400,800,1600,3200\}$. The figure displays both the sample mean and the 95\%-confidence interval calculated by the results of 100 repeats for the regret. Observe that,  even when the problem \eqref{prob:fluid} is degenerate, the regrets are almost constants for full information / informed price settings and $O(\sqrt{T})$ regret in no information case. This matches our theoretical guarantee.

\begin{figure}[h!]
\vskip 0.2in
    \begin{center}
\centerline{\includegraphics[width=0.8\linewidth]{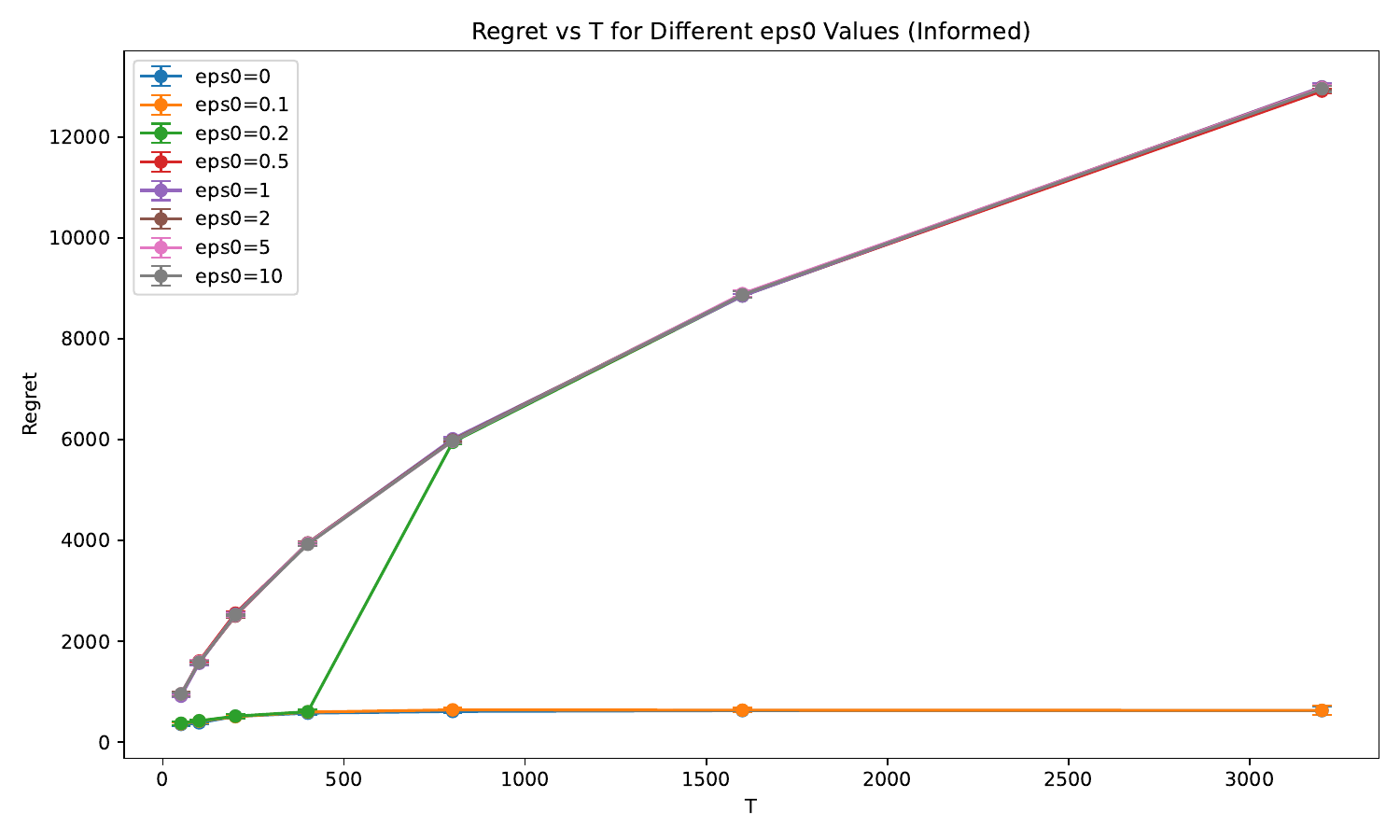}}
\caption{Regret of our algorithms under different number of time horizon $T$ and misspecification error $\epsilon^0$\label{fig:comparison4}.}
\end{center}
\vskip -0.2in
\end{figure}

Figure \ref{fig:comparison4} illustrates the relationship between the regret and the misspecification error $\epsilon^0$ in Algorithm \ref{alg:resolve_learn_incumbent} with $\rho=0.1$. We observe a phase transition as the $\epsilon^0$ grows and the regret exchanges from $O((\epsilon^0)^2T)$ to $O(\sqrt{T})$.

\section{Conclusions}

In this paper, we investigate the dynamic pricing problem with resource constraints by developing algorithms that effectively balance exploration and exploitation. Our main contributions are as follows:

\begin{enumerate}
    \item \textbf{Boundary Attracted Re-solve Method}: We introduce the Boundary Attracted Re-solve Method (Algorithm \ref{alg:resolve}) for the full information setting. This method achieves a logarithmic regret without relying on the non-degeneracy condition typically assumed in previous studies. 
    \item \textbf{Learning Without Prior Information}: Extending our framework to settings with no prior information, we develop an online learning algorithm (Algorithm \ref{alg:resolve_learn}) that attains a regret of $O(\sqrt{T})$. This matches the known lower bounds, demonstrating that our approach is optimal even when the demand function is unknown initially.
    \item \textbf{Incorporating Machine-learned Informed Prices}: We consider the practical availability of offline data and introduce an estimate-then-select re-solve algorithm (Algorithm \ref{alg:resolve_learn_incumbent}) that utilizes informed price-demand pairs with known estimation errors. Our analysis shows that this method bridges the gap between the full information and no information settings, providing improved regret bounds when reliable offline data is available.
\end{enumerate}

Our numerical experiments support the theoretical results, showing that our algorithms perform well across different informational settings. Specifically, our methods maintain low regret even in cases where the problem is degenerate, highlighting their robustness and practical applicability.

Future research could explore extensions to more complex demand models beyond the linear assumption and investigate adaptive techniques for situations where estimation error bounds are uncertain or vary over time. Additionally, integrating advanced machine learning methods for parameter estimation may further enhance the performance of our algorithms in real-world applications.


\bibliographystyle{apalike}

\bibliography{main}

\begin{thebibliography}{}

\bibitem[Abbasi-Yadkori et~al., 2011]{abbasi2011improved}
Abbasi-Yadkori, Y., P{\'a}l, D., and Szepesv{\'a}ri, C. (2011).
\newblock Improved algorithms for linear stochastic bandits.
\newblock {\em Advances in neural information processing systems}, 24.

\bibitem[Agrawal and Devanur, 2016]{agrawal2016linear}
Agrawal, S. and Devanur, N. (2016).
\newblock Linear contextual bandits with knapsacks.
\newblock {\em Advances in neural information processing systems}, 29.

\bibitem[Agrawal et~al., 2016]{agrawal2016efficient}
Agrawal, S., Devanur, N.~R., and Li, L. (2016).
\newblock An efficient algorithm for contextual bandits with knapsacks, and an extension to concave objectives.
\newblock In {\em Conference on Learning Theory}, pages 4--18. PMLR.

\bibitem[Amin et~al., 2014]{amin2014repeated}
Amin, K., Rostamizadeh, A., and Syed, U. (2014).
\newblock Repeated contextual auctions with strategic buyers.
\newblock {\em Advances in Neural Information Processing Systems}, 27.

\bibitem[Ao et~al., 2024a]{ao2024online}
Ao, R., Chen, H., Simchi-Levi, D., and Zhu, F. (2024a).
\newblock Online local false discovery rate control: A resource allocation approach.
\newblock {\em arXiv preprint arXiv:2402.11425}.

\bibitem[Ao et~al., 2024b]{ao2024two}
Ao, R., Fu, H., and Simchi-Levi, D. (2024b).
\newblock Two-stage online reusable resource allocation: Reservation, overbooking and confirmation call.
\newblock {\em arXiv preprint arXiv:2410.15245}.

\bibitem[Auer, 2002]{auer2002using}
Auer, P. (2002).
\newblock Using confidence bounds for exploitation-exploration trade-offs.
\newblock {\em Journal of Machine Learning Research}, 3(Nov):397--422.

\bibitem[Badanidiyuru et~al., 2014]{badanidiyuru2014resourceful}
Badanidiyuru, A., Langford, J., and Slivkins, A. (2014).
\newblock Resourceful contextual bandits.
\newblock In {\em Conference on Learning Theory}, pages 1109--1134. PMLR.

\bibitem[Balseiro et~al., 2020]{balseiro2020dual}
Balseiro, S., Lu, H., and Mirrokni, V. (2020).
\newblock Dual mirror descent for online allocation problems.
\newblock In {\em International Conference on Machine Learning}, pages 613--628. PMLR.

\bibitem[Banerjee and Freund, 2020]{banerjee2020uniform}
Banerjee, S. and Freund, D. (2020).
\newblock sc.
\newblock In {\em Abstracts of the 2020 SIGMETRICS/Performance Joint International Conference on Measurement and Modeling of Computer Systems}, pages 1--2.

\bibitem[Besbes and Zeevi, 2015]{besbes2015surprising}
Besbes, O. and Zeevi, A. (2015).
\newblock On the (surprising) sufficiency of linear models for dynamic pricing with demand learning.
\newblock {\em Management Science}, 61(4):723--739.

\bibitem[Bonnans and Shapiro, 2013]{bonnans2013perturbation}
Bonnans, J.~F. and Shapiro, A. (2013).
\newblock {\em Perturbation analysis of optimization problems}.
\newblock Springer Science \& Business Media.

\bibitem[Boyd and Vandenberghe, 2004]{boyd2004convex}
Boyd, S. and Vandenberghe, L. (2004).
\newblock {\em Convex optimization}.
\newblock Cambridge university press.

\bibitem[Broder and Rusmevichientong, 2012]{broder2012dynamic}
Broder, J. and Rusmevichientong, P. (2012).
\newblock Dynamic pricing under a general parametric choice model.
\newblock {\em Operations Research}, 60(4):965--980.

\bibitem[Bu et~al., 2022]{bu2022context}
Bu, J., Simchi-Levi, D., and Wang, C. (2022).
\newblock Context-based dynamic pricing with partially linear demand model.
\newblock {\em Advances in Neural Information Processing Systems}, 35:23780--23791.

\bibitem[Bu et~al., 2020]{bu2020online}
Bu, J., Simchi-Levi, D., and Xu, Y. (2020).
\newblock Online pricing with offline data: Phase transition and inverse square law.
\newblock In {\em International Conference on Machine Learning}, pages 1202--1210. PMLR.

\bibitem[Bumpensanti and Wang, 2020]{bumpensanti2020re}
Bumpensanti, P. and Wang, H. (2020).
\newblock A re-solving heuristic with uniformly bounded loss for network revenue management.
\newblock {\em Management Science}, 66(7):2993--3009.

\bibitem[Chen et~al., 2021]{chen2021joint}
Chen, Q., Jasin, S., and Duenyas, I. (2021).
\newblock Joint learning and optimization of multi-product pricing with finite resource capacity and unknown demand parameters.
\newblock {\em Operations Research}, 69(2):560--573.

\bibitem[Chen et~al., 2022]{Chen2022}
Chen, Q.~G., Wang, H., and Wang, Z. (2022).
\newblock {\em Learning and Pricing with Inventory Constraints}, volume~18, pages 103--135.
\newblock Springer Nature.

\bibitem[Chen et~al., 2024]{chen2024contextual}
Chen, Z., Ai, R., Yang, M., Pan, Y., Wang, C., and Deng, X. (2024).
\newblock Contextual decision-making with knapsacks beyond the worst case.
\newblock In {\em The Thirty-eighth Annual Conference on Neural Information Processing Systems}.

\bibitem[Cheung and Lyu, 2024]{cheung2024leveraging}
Cheung, W.~C. and Lyu, L. (2024).
\newblock Leveraging (biased) information: Multi-armed bandits with offline data.
\newblock {\em arXiv preprint arXiv:2405.02594}.

\bibitem[Cohen et~al., 2020]{cohen2020feature}
Cohen, M.~C., Lobel, I., and Paes~Leme, R. (2020).
\newblock Feature-based dynamic pricing.
\newblock {\em Management Science}, 66(11):4921--4943.

\bibitem[Dudik et~al., 2011]{dudik2011efficient}
Dudik, M., Hsu, D., Kale, S., Karampatziakis, N., Langford, J., Reyzin, L., and Zhang, T. (2011).
\newblock Efficient optimal learning for contextual bandits.
\newblock {\em arXiv preprint arXiv:1106.2369}.

\bibitem[Elmachtoub and Grigas, 2022]{elmachtoub2022smart}
Elmachtoub, A.~N. and Grigas, P. (2022).
\newblock Smart “predict, then optimize”.
\newblock {\em Management Science}, 68(1):9--26.

\bibitem[Ferreira et~al., 2018]{ferreira2018online}
Ferreira, K.~J., Simchi-Levi, D., and Wang, H. (2018).
\newblock Online network revenue management using thompson sampling.
\newblock {\em Operations research}, 66(6):1586--1602.

\bibitem[Gallego and Van~Ryzin, 1994]{gallego1994optimal}
Gallego, G. and Van~Ryzin, G. (1994).
\newblock Optimal dynamic pricing of inventories with stochastic demand over finite horizons.
\newblock {\em Management science}, 40(8):999--1020.

\bibitem[Grossman and Stiglitz, 1976]{grossman1976information}
Grossman, S.~J. and Stiglitz, J.~E. (1976).
\newblock Information and competitive price systems.
\newblock {\em The American economic review}, 66(2):246--253.

\bibitem[Jaillet et~al., 2024]{jaillet2024should}
Jaillet, P., Podimata, C., Vakhutinsky, A., and Zhou, Z. (2024).
\newblock When should you offer an upgrade: Online upgrading mechanisms for resource allocation.
\newblock {\em arXiv preprint arXiv:2402.08804}.

\bibitem[Jasin, 2014]{jasin2014reoptimization}
Jasin, S. (2014).
\newblock Reoptimization and self-adjusting price control for network revenue management.
\newblock {\em Operations Research}, 62(5):1168--1178.

\bibitem[Jasin and Kumar, 2012]{jasin2012re}
Jasin, S. and Kumar, S. (2012).
\newblock A re-solving heuristic with bounded revenue loss for network revenue management with customer choice.
\newblock {\em Mathematics of Operations Research}, 37(2):313--345.

\bibitem[Jasin and Kumar, 2013]{jasin2013analysis}
Jasin, S. and Kumar, S. (2013).
\newblock Analysis of deterministic lp-based booking limit and bid price controls for revenue management.
\newblock {\em Operations Research}, 61(6):1312--1320.

\bibitem[Javanmard and Nazerzadeh, 2019]{javanmard2019dynamic}
Javanmard, A. and Nazerzadeh, H. (2019).
\newblock Dynamic pricing in high-dimensions.
\newblock {\em Journal of Machine Learning Research}, 20(9):1--49.

\bibitem[Jiang et~al., 2022]{jiang2022degeneracy}
Jiang, J., Ma, W., and Zhang, J. (2022).
\newblock Degeneracy is ok: Logarithmic regret for network revenue management with indiscrete distributions.
\newblock {\em arXiv preprint arXiv:2210.07996}.

\bibitem[Keskin and Zeevi, 2014]{keskin2014dynamic}
Keskin, N.~B. and Zeevi, A. (2014).
\newblock Dynamic pricing with an unknown demand model: Asymptotically optimal semi-myopic policies.
\newblock {\em Operations research}, 62(5):1142--1167.

\bibitem[Kumar and Kleinberg, 2022]{kumar2022non}
Kumar, R. and Kleinberg, R. (2022).
\newblock Non-monotonic resource utilization in the bandits with knapsacks problem.
\newblock {\em Advances in Neural Information Processing Systems}, 35:19248--19259.

\bibitem[Lattimore and Szepesv{\'a}ri, 2020]{lattimore2020bandit}
Lattimore, T. and Szepesv{\'a}ri, C. (2020).
\newblock {\em Bandit algorithms}.
\newblock Cambridge University Press.

\bibitem[Li and Ye, 2022]{li2022online}
Li, X. and Ye, Y. (2022).
\newblock Online linear programming: Dual convergence, new algorithms, and regret bounds.
\newblock {\em Operations Research}, 70(5):2948--2966.

\bibitem[Li et~al., 2021]{li2021unifying}
Li, Y., Xie, H., Lin, Y., and Lui, J.~C. (2021).
\newblock Unifying offline causal inference and online bandit learning for data driven decision.
\newblock In {\em Proceedings of the Web Conference 2021}, pages 2291--2303.

\bibitem[Liu and Grigas, 2022]{liu2022online}
Liu, H. and Grigas, P. (2022).
\newblock Online contextual decision-making with a smart predict-then-optimize method.
\newblock {\em arXiv preprint arXiv:2206.07316}.

\bibitem[Liu et~al., 2022]{liu2022non}
Liu, S., Jiang, J., and Li, X. (2022).
\newblock Non-stationary bandits with knapsacks.
\newblock {\em Advances in Neural Information Processing Systems}, 35:16522--16532.

\bibitem[Reiman and Wang, 2008]{reiman2008asymptotically}
Reiman, M.~I. and Wang, Q. (2008).
\newblock An asymptotically optimal policy for a quantity-based network revenue management problem.
\newblock {\em Mathematics of Operations Research}, 33(2):257--282.

\bibitem[Sankararaman and Slivkins, 2021]{sankararaman2021bandits}
Sankararaman, K.~A. and Slivkins, A. (2021).
\newblock Bandits with knapsacks beyond the worst case.
\newblock {\em Advances in Neural Information Processing Systems}, 34:23191--23204.

\bibitem[Shah et~al., 2019]{shah2019semi}
Shah, V., Johari, R., and Blanchet, J. (2019).
\newblock Semi-parametric dynamic contextual pricing.
\newblock {\em Advances in Neural Information Processing Systems}, 32.

\bibitem[Sivakumar et~al., 2020]{sivakumar2020structured}
Sivakumar, V., Wu, S., and Banerjee, A. (2020).
\newblock Structured linear contextual bandits: A sharp and geometric smoothed analysis.
\newblock In {\em International Conference on Machine Learning}, pages 9026--9035. PMLR.

\bibitem[Sivakumar et~al., 2022]{sivakumar2022smoothed}
Sivakumar, V., Zuo, S., and Banerjee, A. (2022).
\newblock Smoothed adversarial linear contextual bandits with knapsacks.
\newblock In {\em International Conference on Machine Learning}, pages 20253--20277. PMLR.

\bibitem[Vera and Banerjee, 2021]{vera2021bayesian}
Vera, A. and Banerjee, S. (2021).
\newblock The bayesian prophet: A low-regret framework for online decision making.
\newblock {\em Management Science}, 67(3):1368--1391.

\bibitem[Wainwright, 2019]{wainwright2019high}
Wainwright, M.~J. (2019).
\newblock {\em High-dimensional statistics: A non-asymptotic viewpoint}, volume~48.
\newblock Cambridge university press.

\bibitem[Wang and Wang, 2022]{wang2022constant}
Wang, Y. and Wang, H. (2022).
\newblock Constant regret resolving heuristics for price-based revenue management.
\newblock {\em Operations Research}, 70(6):3538--3557.

\bibitem[Wang and Zheng, 2023]{wang2023measuring}
Wang, Y. and Zheng, Z. (2023).
\newblock Measuring policy performance in online pricing with offline data: Worst-case perspective and bayesian perspective.
\newblock {\em Journal of Systems Science and Systems Engineering}, 32(3):352--371.

\bibitem[Wang et~al., 2024]{wang2024online}
Wang, Y., Zheng, Z., and Shen, Z.-J.~M. (2024).
\newblock Online pricing with polluted offline data.
\newblock {\em Available at SSRN 4320324}.

\bibitem[Xu and Wang, 2021]{xu2021logarithmic}
Xu, J. and Wang, Y.-X. (2021).
\newblock Logarithmic regret in feature-based dynamic pricing.
\newblock {\em Advances in Neural Information Processing Systems}, 34:13898--13910.

\bibitem[Xu and Wang, 2022]{xu2022towards}
Xu, J. and Wang, Y.-X. (2022).
\newblock Towards agnostic feature-based dynamic pricing: Linear policies vs linear valuation with unknown noise.
\newblock In {\em International Conference on Artificial Intelligence and Statistics}, pages 9643--9662. PMLR.

\bibitem[Zhang and Cheung, 2024]{zhang2024piecewisestationary}
Zhang, X. and Cheung, W.~C. (2024).
\newblock Piecewise-stationary bandits with knapsacks.
\newblock In {\em The Thirty-eighth Annual Conference on Neural Information Processing Systems}.

\end{thebibliography}

\newpage
\appendix
\section{Proof of Theorem \ref{thm:full_information}\label{appendix:resolve}}
To begin with, we restate the re-solve constrained programming problem:
\begin{equation}
    \label{prob:resolve_restate}
    \begin{aligned}
        \max_{p\in\Pcal}\qquad & r ={\bm{p}}^\top{\bm{d}}\\ 
    \st\qquad& {\bm{d}} = {\bm{\alpha}} + B{\bm{p}},\\ 
    &Ad\le \frac{{\bm{c}}^t}{T-t+1},
    \end{aligned}
\end{equation}
where ${\bm{c}}^t$ is the inventory level at the beginning of time $t$. We denote $(p^{\pi,t},{\bm{d}}^{\pi,t})$ as the optimal solution to \eqref{prob:resolve_restate}. For notational convenience, we use $p({\bm{d}}):=B^{-1}({\bm{d}}-{\bm{\alpha}}),r({\bm{d}}):=p({\bm{d}})^\top{\bm{d}}$ to represent the unnoised maps of demand to price and demand to revenue, respectively. Moreover, we use $r({\bm{d}},{\bm{\epsilon}}):=r({\bm{d}})+p({\bm{d}})^\top{\bm{\epsilon}}$ to represent the revenue with noised demand. Throughout the proof, we use $\pi$ to denote our online policy given in Algorithm \ref{alg:resolve}. Given filtration $\F^T$ and policy $\pi'$, we use $\Rcal^T(\pi',\F^T)$ to represent the total revenue under policy $\pi'$ and realized sample path $\F^T$. Now we introduce the following concept of \textbf{Hybrid} policy, which is crucial in our proof. 
\begin{definition}
    For $1\le t \le T+1$, we define $\mix^{t}$ as the policy that applies online policy $\pi$ in time $1,\dots,t-1$ and ${\bm{d}}^{\pi,t}$, while for the periods $[t+1,T]$, there are no noises of demand. Moreover, define $\mix^{1}$ as the fluid optimal policy given in \eqref{prob:fluid} without noises and $\mix^{T+1} = \pi$ as the online policy $\pi$ throughout the process.
\end{definition}
By definition, we have 
\begin{equation*}
    \Rcal^T(\mix^{t},\F^T) = \sum_{s=1}^{t-1} r({\bm{d}}^{\pi,s},{\bm{\epsilon}}^s) + (T-t+1)r({\bm{d}}^{\pi,t}).
\end{equation*}
Therefore, it holds that $\ex{}{\Rcal^T(\mix^{t},\F^T)}\ge \ex{}{f(\mix^{t+1},\F^T)},0\le t\le T-1$ by the convexity of \eqref{prob:resolve_restate}. The regret can then be decomposed as follows:
\begin{align}
    \label{eq:regret_decomp1}
    \regret[T]{\pi} &= \ex{}{\sum_{t=1}^{T}\Rcal^T(\mix^{t},\F^T)-\Rcal^T(\mix^{t+1},\F^T)}\nonumber\\
    &=\sum_{t=1}^{T}\ex{}{\Rcal^T(\mix^{t},\F^T)-\Rcal^T(\mix^{t+1},\F^T)}.
\end{align}
We now focus on how to give bound to the term $\ex{}{\Rcal^T(\mix^{t},\F^T)-\Rcal^T(\mix^{t+1},\F^T)}$ for $0\le t\le T-1$. 

By definition, the realized demand at time $t$ is given by ${\bm{d}}^{t} = {\bm{d}}^{\pi,t} + {\bm{\epsilon}}^t$ and we have 
\begin{equation}
    \label{eq:taylor_expansion}
    \begin{aligned}
    r({\bm{d}}^{\pi,t},{\bm{\epsilon}}^t) &= r({\bm{d}}^{\pi,t}) +({\bm{\epsilon}}^t)^\top{\bm{p}}^{\pi,t},\\
    r(d') &= r({\bm{d}})+(d'-d)^\top B^{-1}(d'-d)+(2d-{\bm{\alpha}})^\top B^{-1}(d'-d),\quad\forall {\bm{d}},d'\in\R_+^n
    \end{aligned}
\end{equation}
 The single-step difference can then be rewritten as:
\begin{equation}
    \label{eq:single_step1}
    \begin{aligned}
        \Rcal^T(\mix^{t},\F^T)-\Rcal^T(\mix^{t+1},\F^T) &= (T-t+1)r({\bm{d}}^{\pi,t}) - r({\bm{d}}^{\pi,t},{\bm{\epsilon}}^t) - (T-t)r({\bm{d}}^{\pi,t+1})\\ 
    \end{aligned}
\end{equation}
Now we proceed with the following three cases: case (I) $\min_id_i^{\pi,t}\ge\zeta(T-t+1)^{-1/2}$; case (II) $\max_i{\bm{d}}^{\pi,t} \le \zeta(T-t+1)^{-1/2}$; (III) $\min_id_i^{\pi,t}\le \zeta(T-t+1)^{-1/2}\le \max_id_i^{\pi,t}.$

\paragraph*{Case (I) $\min_id_i^{\pi,t}> \zeta(T-t+1)^{-1/2} $.} For convenience, we omit the notation for conditional expectation / probability of the event $\{{\bm{d}}^{\pi,t}> \zeta(T-t+1)^{-1/2}\}$. 
Denote $\Ecal_i^t $ as the event $d_i^{\pi,t} \ge \frac{\epsilon_i^t}{T-t},\forall i\in[n]$ and $\Ecal^t=\cap_{i=1}^n\Ecal_i^t$. We now split the \eqref{eq:single_step1} according to 
\begin{equation}\label{eq:decomp10}
    \begin{aligned}
        &\ex{}{\Rcal^{T}(\mix^{t},\F^T)-\Rcal^T(\mix^{t+1},\F^T)\big\vert\F^{t-1}} 
    \\&= \P(\Ecal^t)\ex{}{\Rcal^T(\mix^{t},\F^T)-\Rcal^T(\mix^{t+1},\F^T)\big\vert\Ecal^t,\F^{t}}\\&\qquad+ \P((\Ecal^t)^c)\ex{}{\Rcal^T(\mix^{t},\F^T)-\Rcal^T(\mix^{t+1},\F^T)\big\vert(\Ecal^t)^c,\F^{t-1}}.
    \end{aligned}
\end{equation}
For the first term in \eqref{eq:decomp10}, condition on $\Ecal^t$, by definition we have 
\begin{equation}
    \label{eq:decomp11}
    \begin{aligned}
        &\ex{}{\Rcal^T(\mix^{t},\F^T)-\Rcal^T(\mix^{t+1},\F^T)\big\vert\Ecal^t,\F^{t-1}} \\ 
        &= \ex{}{(T-t+1)r({\bm{d}}^{\pi,t})-r({\bm{d}}^{\pi,t},{\bm{\epsilon}}^t)   - (T-t)r({\bm{d}}^{\pi,t+1})\big\vert\Ecal^t,\F^{t-1}}\\ 
        &\le \ex{}{(T-t+1)r({\bm{d}}^{\pi,t})-r({\bm{d}}^{\pi,t},{\bm{\epsilon}}^t) - (T-t)r({\bm{d}}^{\pi,t}-\frac{{\bm{\epsilon}}^t}{T-t})\big\vert\Ecal^t,\F^{t-1}}\\ 
        &\overset{(a)}{=}\ex{}{-({\bm{\epsilon}}^t)^\top{\bm{p}}^{\pi,t}+ (T-t)\frac{({\bm{\epsilon}}^t)^\top}{T-t}B^{-1}\frac{{\bm{\epsilon}}^t}{T-t}+(T-t)(2{\bm{d}}^{\pi,t}-{\bm{\alpha}})^\top B^{-1}\frac{{\bm{\epsilon}}^t}{T-t} |\Ecal^t,\F^{t-1}}\\
        &\overset{(b)}{=} \ex{}{({\bm{d}}^{\pi,t})^\top B^{-1}{\bm{\epsilon}}^t+\frac{1}{T-t}({\bm{\epsilon}}^t)^\top B^{-1}{\bm{\epsilon}}^t\vert \Ecal^t,\F^{t-1}},
    \end{aligned}
\end{equation}
where $(a)$ follows from a similar argument in \eqref{eq:taylor_expansion} and $(b)$ follows from the definition: $p^{\pi,t}=B^{-1}({\bm{d}}^{\pi,t}-{\bm{\alpha}})$. We then have 
\begin{equation}
    \label{eq:decomp12}
    \begin{aligned}
        &\P(\Ecal^t)\ex{}{({\bm{d}}^{\pi,t})^\top B^{-1}{\bm{\epsilon}}^t|\Ecal^t,\F^{t-1}}\\ 
        & = \ex{}{({\bm{d}}^{\pi,t})^\top B^{-1}{\bm{\epsilon}}^t}-\P((\Ecal^t)^c)\ex{}{({\bm{d}}^{\pi,t})^\top B^{-1}{\bm{\epsilon}}^t|(\Ecal^t)^c,\F^{t-1}}\\ 
        &= -\ex{}{({\bm{d}}^{\pi,t})^\top B^{-1}{\bm{\epsilon}}^t\mathbbm 1\{(\Ecal)^c\}|\F^{t-1}}\\ 
        &\overset{(a)}{\le} \norm{B^{-1}}_2d_{\max}\P((\Ecal)^c\vert\F^{t-1})^{1/2}\ex{}{\norm{{\bm{\epsilon}}^t}_2
        ^2\vert\F^{t-1}}^2\\ 
        &= \norm{B^{-1}}_2d_{\max}\P((\Ecal)^c\vert\F^{t-1})^{1/2}\ex{}{\norm{{\bm{\epsilon}}^t}_2^2}^{1/2}
    \end{aligned}
\end{equation}
where we applied Cauchy-Schwarz inequality in (a). Similarly, we have 
\begin{equation}\label{eq:decomp13}
    \begin{aligned}
        \P(\Ecal^t)\ex{}{-\frac{({\bm{\epsilon}}^t)^\top B^{-1}{\bm{\epsilon}}^t}{T-t}\big\vert\Ecal^t,\F^{t-1}} & \le  \frac{\norm{B^{-1}}_2}{T-t}\ex{}{({\bm{\epsilon}}^t)^2}.
    \end{aligned}
\end{equation}

On the other hand, with a similar argument, we can give bound to the term 
\begin{equation}
    \label{eq:decomp14}
    \begin{aligned}
        &\P((\Ecal^t)^c)\ex{}{f(\mix^{t-1})-f(\mix^{t})\big\vert(\Ecal^t)^c} \\ 
        &\le \P((\Ecal^t)^c)\ex{}{(T-t)r({\bm{d}}^{\pi,d})t)r({\bm{d}}^{\pi,t})-({\bm{\epsilon}}^t)^\top{\bm{p}}^{\pi,t}\vert (\Ecal^t)^c}\\
        &\le \P((\Ecal^t)^c)\prn{(T-t)r_{\max}-\ex{}{({\bm{\epsilon}}^t)^\top{\bm{p}}^{\pi,t}\vert(\Ecal^t)^c}}\\
        &\le (T-t)r_{\max}\P((\Ecal^t)^c)+\sqrt{n}U\ex{}{({\bm{\epsilon}}^t)^2}^{1/2}\P((\Ecal^t)^c)^{1/2},
    \end{aligned}
\end{equation}
where we again use Cauchy-Schwarz inequality in the last line.

Now combining \eqref{eq:decomp10} \eqref{eq:decomp11} \eqref{eq:decomp12} \eqref{eq:decomp13} \eqref{eq:decomp14}, we get 
\begin{align}
    &\label{eq:single_step2}
    \ex{}{\Rcal^T(\mix^t,\F^T)-\Rcal^T(\mix^{t+1},\F^T)}\nonumber\\
    & \le\sigma(\norm{B^{-1}}_2d_{\max}+\sqrt{n}U)\P((\Ecal^t)^c)^{1/2}+(T-t)r_{\max}\P((\Ecal^t)^c)+\frac{\sigma^2\norm{B^{-1}}_2}{T-t},
\end{align}
Now in order to give bound to $\P((\Ecal^t)^c),$ we introduce the following concentration inequality.
\begin{lemma}[\citealt{wainwright2019high}]\label{lem:tail}
    Let $X_1,\dots,X_n$ be $\sigma^2$-sub-Gaussian random variables with zero mean, then for each $\lambda >0,$ it holds that 
    \begin{align*}
        \P\prn{\max_{1\le i\le n}X_i\ge \lambda}\le n\exp(-\lambda/2\sigma^2)
    \end{align*}
\end{lemma}
By Lemma \ref{lem:tail}, we can bound the term $\P((\Ecal^t)^c)$ as
\begin{align}
\label{eq:prob_complementary}
        \P((\Ecal^t)^c) &= \P(\exists i\in[n], \st d_i^{\pi,[t,T]} < \frac{\epsilon_i^t}{T-t})\nonumber\\ 
        &\le \P(\frac{\zeta}{(T-t+1)^{1/2}}<\max_i\frac{\epsilon_i^t}{T-t})\nonumber\\
        &\overset{(a)}{\le} n\exp(-2(T-t)\log n)\nonumber\\ 
        & = \exp(-2(T-t)),
\end{align}
where (a) is detived from Lemma \ref{lem:tail} and the fact that $\min_id_i^{\pi,t}\ge\zeta(T-t+1)^{-1/2} $ and $\zeta\ge 4\sigma^2\log n$. Plugging the above inequality into \eqref{eq:single_step2} leads to 
\begin{equation}\label{eq:eq:single_step21}
\begin{aligned}
    &\ex{}{\Rcal^T(\mix^{t},\F^{T})-\Rcal^T(\mix^{t+1},\F^T)}\\ 
    &\le \sigma(\norm{B^{-1}}_2d_{\max}+\sqrt{n}U)\exp(-(T-t))+r_{\max}(T-t)\exp(-2(T-t))+\frac{\sigma^2\norm{B^{-1}}_2}{T-t}.
\end{aligned}
\end{equation}

\paragraph*{Case (II): $\max_id_i^{\pi,}\le \zeta/(T-t+1)^{-1/2}$.} In this case, we have ${\bm{d}}^{\pi,[t,T]} = {\bm{\epsilon}}^t = 0$. Using \eqref{eq:single_step1}, the single-step difference in case (II) can be upper bounded by 
\begin{equation}
    \label{eq:single_step22}
    \begin{aligned}
        &\ex{}{\Rcal^T(\mix^{t},\F^T)-\Rcal^{T}(\mix^{t+1},\F^T)}\\ 
                &\le \ex{}{(T-t+1)r({\bm{d}}^{\pi,t})-(T-t)r\prn{\frac{T-t+1}{T-t}{\bm{d}}^{\pi,t}}  }\\ 
        &= \ex{}{(T-t+1)({\bm{d}}^{\pi,t})^\top B^{-1}({\bm{d}}^{\pi,t}-{\bm{\alpha}})- (T-t)\frac{T-t+1}{T-t}({\bm{d}}^{\pi,t})^\top B^{-1}(\frac{T-t+1}{T-t}{\bm{d}}^{\pi,t}-{\bm{\alpha}} )   }\\ 
        &= \ex{}{-\frac{T-t+1}{T-t}({\bm{d}}^{\pi,t})^\top B^{-1}{\bm{d}}^{\pi,t}}\\ 
        &\le \frac{2n^2\zeta^2\norm{B^{-1}}_2}{T-t+1}
    \end{aligned}
\end{equation}

\paragraph*{Case (III): $\min_id_i^{\pi,t}<\zeta(T-t+1)^{-1/2}<\max_id_i^{\pi,t}.$} In this case, we can derive upper bound by following both Case (I) and Case (II). To be more specific, let $\Ical = \{i\in[n]:d_i^{\pi,t}>\zeta(T-t+1)^{-1/2}\}$, $\overline{\Ical}=[n]\backslash\Ical$ and $\Ecal_\Ical = \cap_{i\in \Ical}\Ecal_i.$ Let $\epsilon_\Ical^t$ be the vector that has components $\epsilon_i^t$ for $i\in\Ical$ and $0$ otherwise. We let $\tilde {\bm{d}}^t$ be vector with components $d_i^{\pi,t}-{\bm{\epsilon}}^t/(T-t),i\in\Ical$ and $(T-t+1)d_i^{\pi,t}/(T-t),i\in\overline\Ical.$ Then the single-step difference can be upper bounded by
\begin{align}
\label{eq:decomp_15}
    &\ex{}{\Rcal^T(\mix^t,\F^T-\Rcal^T(\mix^{t+1},\F^T)}\nonumber\\ 
    &\le \ex{}{(T-t+1)r({\bm{d}}^{\pi,t})-r({\bm{d}}^t,{\bm{\epsilon}}^t)-(T-t)r(\tilde {\bm{d}}^t)}\nonumber\\ 
    & = \ex{}{(p^t)^\top{\bm{\epsilon}}^t-({\bm{d}}^t-{\bm{d}}^{\pi,t})^\top B^{-1}({\bm{d}}^t-{\bm{d}}^{\pi,t})-(2{\bm{d}}^{t,\pi}-{\bm{\alpha}})^\top B^{-1}({\bm{d}}^t-{\bm{d}}^{\pi,t}) }\nonumber\\
    &\qquad - (T-t)\ex{}{(\tilde {\bm{d}}^t-{\bm{d}}^{\pi,t})^\top B^{-1}(\tilde {\bm{d}}^t-{\bm{d}}^{\pi,t})+(2{\bm{d}}^{\pi,t}-{\bm{\alpha}})^\top B^{-1}(\tilde {\bm{d}}^t-{\bm{d}}^{\pi,t})}.
\end{align}
Note that $d_i^t = d_i^{\pi,t}$ for $i\in\Ical$ and $d_i^t = 0$ for $i\in\overline\Ical$. Moreover, $d_i^{\pi,t}\le \zeta(T-t+1)^{-1/2},i\in\overline\Ical.$ We have 
\begin{align}\label{eq:case_3_1}
    ({\bm{d}}^t-{\bm{d}}^{\pi,t})^\top B^{-1}({\bm{d}}^t-{\bm{d}}^{\pi,t})\le \zeta^2(n-|\Ical|)\norm{B^{-1}}_2(T-t+1)^{-1}.
\end{align}
On the other hand, we have 
\begin{align}
    \label{eq:case_3_2}
    &(\tilde {\bm{d}}^t-{\bm{d}}^{\pi,t})^\top B^{-1}(\tilde {\bm{d}}^t-{\bm{d}}^{\pi,t})\nonumber\\ 
    &= \frac{1}{(T-t)^2}({\bm{d}}_{\overline\Ical}^{\pi,t})^\top B^{-1}{\bm{d}}_{\overline\Ical}^{\pi,t}-\frac{2}{(T-t)^2}({\bm{d}}_{\overline\Ical}^{\pi,t})^\top B^{-1}\epsilon_\Ical^t+\frac{1}{(T-t)^2}(\epsilon_\Ical^t)^\top B^{-1}\epsilon_\Ical^t.
\end{align}
Moreover, for the first-order terms, we have 
\begin{align}
    \label{eq:case_3_3}
    (2{\bm{d}}^{t,\pi}-{\bm{\alpha}})^\top B^{-1}({\bm{d}}^t-{\bm{d}}^{\pi,t})=-(2{\bm{d}}^{t,\pi}-{\bm{\alpha}})^\top B^{-1}{\bm{d}}_{\overline{\Ical}}^{\pi,t},
\end{align}
and
\begin{align}
    \label{eq:case_3_4}
    (2{\bm{d}}^{\pi,t}-{\bm{\alpha}})B^{-1}(\tilde {\bm{d}}^t-{\bm{d}}^{\pi,t}) = \frac{1}{T-t}(2{\bm{d}}^{\pi,t}-{\bm{\alpha}})^\top B^{-1}{\bm{d}}^{\pi,t} - \frac{1}{T-t}(2{\bm{d}}^{\pi,t}-{\bm{\alpha}})^\top B^{-1}\epsilon_\Ical^\top.
\end{align}
Plugging \eqref{eq:case_3_1}, \eqref{eq:case_3_2}, \eqref{eq:case_3_3} and \eqref{eq:case_3_4} into \eqref{eq:decomp_15} yields as similar bound as in Case (I) and Case (II):
\begin{align}
    \label{eq:single_step23}
    &\ex{}{\Rcal^T(\mix^t,\F^T)-\Rcal^T(\mix^{t+1},\F^T)}\nonumber\\ 
    &\le \sigma(\norm{B^{-1}}_2d_{\max}+\sqrt{n}U)\exp(-(T-t))+r_{\max}(T-t)\exp(-2(T-t)) + \frac{2n^2\zeta^2\norm{B^{-1}}_2}{T-t+1}.
\end{align}

\paragraph*{Wrap-up.}
Combining \eqref{eq:eq:single_step21}, \eqref{eq:single_step22} and \eqref{eq:single_step23} leads to 
\begin{align*}
    &\ex{}{\Rcal^T(\mix^{t},\F^T)-\Rcal^T(\mix^{t+1},\F^T}\\ 
    &\le \frac{2n^2\zeta^2\norm{B^{-1}}_2}{T-t+1} + C_0\exp(-(T-t)) + \frac{\sigma^2\norm{B^{-1}}_2}{T-t},
\end{align*}
where $C_0 = C'\sigma(\norm{B^{-1}}_2d_{\max}+\sqrt nU)+r_{\max}$ for some absolute constant $C'$. Plugging them into \eqref{eq:regret_decomp1} yields
\begin{align*}
        \regret[T]{\pi} &=\sum_{t=1}^T\ex{}{\Rcal^T(\mix^{t},\F^T)-\Rcal^T(\mix^{t+1},\F^T)} \\ 
        &=O(2\zeta^2n^2\norm{B^{-1}}_2\log(T) + C_0).
\end{align*}

\section{Proof of Theorem \ref{alg:resolve_learn}}\label{appendix:learn}
In this section, we follow a similar streamline in Appendix \ref{appendix:resolve}. To begin with, recall our re-solve constrained programming problem:
\begin{equation}
    \label{prob:resolve_restate2}
    \begin{aligned}
        \max_{p\in\Pcal}\qquad & r ={\bm{p}}^\top{\bm{d}}\\ 
    \st\qquad& {\bm{d}} = {\bm{\alpha}} + B{\bm{p}},\\ 
    &Ad\le \frac{{\bm{c}}^t}{T-t+1},
    \end{aligned}
\end{equation}
where ${\bm{c}}^t$ is the inventory level at the beginning of time $t$. For time $t$ and $k = \lfloor (t-1)/n\rfloor,$ we use the linear regression \eqref{eq:linear_regression} to fit the coefficients $\hat{\bm{\alpha}}^{kn+1},\hat B^{kn+1}$ and then substitute them into \eqref{prob:resolve_restate2} with inventory level ${\bm{c}}^{kn+1}$ to calculate price $\hat{\bm{p}}^k.$ WLOG we assume that $T=T'n$ for some integer $T'$. Since we're now using a periodic-review re-solve policy, we proceed by modifying the single-step difference in \eqref{eq:regret_decomp1}:
\begin{align}
    \label{eq:regret_decomp2}
    \regret[T]{\pi} &= \ex{}{\sum_{k=1}^{T'}\Rcal^T(\mix^{kn+1},\F^T)-\Rcal^T(\mix^{(k+1)n+1},\F^T)}\nonumber\\&=\sum_{k=1}^{T'}\ex{}{\Rcal^T(\mix^{kn+1},\F^T)-\Rcal^T(\mix^{(k+1)n+1},\F^T)}.
\end{align}
For brevity, we use $k,k+1$ to replace the superscript $kn+1,(k+1)n+1$, respectively when the context is clear. We now focus on how ot give bound to term $\ex{}{\Rcal^T(\mix^{k},\F^T)-\Rcal^T(\mix^{k+1},\F^T)},\forall k\in[T']$.

Recall that for $t\in[T]$ and $k=\lfloor (t-1)/n\rfloor$, $\bm{p}^t, \bm{d}^t:=f(\bm{p}^t)$ are the prices and demands at time $t$ by Algorithm \ref{alg:resolve_learn}. We use $\Delta^t = {\bm{d}}^t- {\bm{d}}^{\pi,k}$ to denote the difference between the targeted demand and the accurate fluid optimal demands by solving \eqref{prob:resolve} with inventory level ${\bm{c}}^{kn+1}.$ Then $\Delta^t$ can be decomposed into two parts of errors as 
\begin{align}
    \label{eq:error_decomp}
    \Delta^t = \underset{:=\Delta_I^t}{\underbrace{(\tilde {\bm{d}}^k-{\bm{d}}^{\pi, nk+1})}}+\underset{:=\Delta_{II}^t}{\underbrace{\overline {\bm{d}}^{t-1}-\overline {\bm{d}}^{kn} }} + \underset{:=\Delta_{III}^t}{\underbrace{\sigma_0t^{-1/4}Be_{t-kn}}}.
\end{align}
Here $\Delta_I^t$ is the estimation error, $\Delta_{II}^t$ is the shift of mean price and $\Delta_{III}^t$ is the perturbation error. With triangular inequality, for $\Delta_{II}^t$ we have
\begin{align}
    \label{eq:error_delta_1}
    \norm{\Delta_{II}^t}_2 &= 
    \norm{\overline{\bm{d}}^{t-1}-\overline{\bm{d}}^{kn}}_2 \nonumber\\
    &\le \frac{t-1-kn}{(t-1)}\norm{\overline{\bm{d}}^{kn}}_2+\frac{1}{t-1}\norm{\sum_{s=kn+1}^{t-1}{\bm{d}}^s}_2\nonumber\\ 
    &\le \frac{n}{t-1}{\bm{d}}_{\max} + \frac{n}{t-1}{\bm{d}}_{\max}\nonumber\\ 
    &=\frac{2n}{t-1}{\bm{d}}_{\max}.
\end{align}
For $\Delta_{III}^t,$ it directly follows from the definition that 
\begin{align}
    \label{eq:error_delta_2}
    \norm{\Delta_{III}^t}_2 &\le \sigma_0\norm{B}_2t^{-1/4}.
\end{align}
Now we proceed by giving bound to the term $\Delta_I^t.$ We define a quantity $J^t$ that is crucial for deriving tail bound for the estimation error. Formally, we let
\begin{align}
    \label{eq:J_t}
    J^k &:= n^{-1}\sum_{s=1}^{nk}(1-s^{-1})\norm{p^s-\overline{\bm{p}}^{s-1}-(\tilde{\bm{p}}^k-\overline{\bm{p}}^{kn})}_2^2\nonumber\\
    &= n^{-1}\sum_{l=1}^k\sum_{i=1}^n\prn{1-\frac{1}{n(l-1)+i}}\sigma_0^2(kn+i)^{-1/2}\nonumber\\ 
    &\ge \frac{\sigma_0^2\sqrt{kn}}{8n},
    \quad\forall t=1,2,\dots.
\end{align}
Now we introduce the following lemma for giving bound to the parameter error:
\begin{lemma}[\citealt{keskin2014dynamic}]
    \label{lem:tail}
    Under the our choice of $p^t,$ there exists constant $C_1,\sigma_1$ such that  
    \begin{align*}
        \P\prn{\norm{\hat{\bm{\alpha}}^{kn+1}-{\bm{\alpha}}}_2+\norm{\hat B^{kn+1}-B}_2>\lambda,J^k\ge \lambda'}\le C_1(kn)^{n^2+n-1}\exp(-\sigma_1(\lambda\wedge\lambda^2)\lambda'),\forall \lambda,\lambda'>0.
    \end{align*}
\end{lemma}
As a result, by combining \eqref{eq:J_t} and Lemma \ref{lem:tail}, we get the following bound:
\begin{align}
    \label{eq:estimation_error_parameter1}
    \P(\norm{\hat{\bm{\alpha}}^{kn+1}-{\bm{\alpha}}}_2+\norm{\hat B^{kn+1}-B}_2>\lambda)\le C_1(kn)^{n^2+n-1}\exp\prn{-\frac{\sigma_0^2\sigma_1\sqrt{kn}}{8n}(\lambda\wedge\lambda^2)}.
\end{align}

Now with a similar argument as in the Proof of Theorem 6 in \citealt{keskin2014dynamic}, we arrive at 
\begin{align}
    \label{eq:estimation_error_parameter2}
    \ex{}{(\norm{\hat{\bm{\alpha}}^{kn+1}-{\bm{\alpha}}}_2+\norm{\hat B^{kn+1}-B}_2)^2}\le \frac{20C_1n^5\log(kn+1)}{\sigma_1\sqrt{kn}}.
\end{align}
In order to use the parameter error to bound the solution error $\Delta_I^t,$ we introduce the
following Lemma concerning the continuity of constrained strongly convex optimization problems.
\begin{lemma}[Prop 4.32, \cite{bonnans2013perturbation}]\label{lem:Lip1}
    Suppose the constraint optimization problem
    \begin{equation*}
        \begin{aligned}
            \max_{p\in\Pcal} \qquad &r({\bm{d}})={\bm{p}}^\top{\bm{d}}\\ 
            \st \qquad & {\bm{d}} = {\bm{\alpha}}+B{\bm{p}},\\
            &Ad \le \frac{{\bm{c}}^t}{T-t+1}\\ 
            &{\bm{d}} \ge 0,
        \end{aligned}
    \end{equation*}
    satisfies second-order growth condition $r({\bm{d}}) \le r({\bm{d}}^{\pi,t}) - \kappa(\dist{{\bm{d}},D^{\pi,t}})^2$ for any $d$ in the feasible set and the optimal solution set $D^{\pi,t}$. Then for any optimal solution $\hat{\bm{d}}$ to the quadratic programming
    \begin{equation*}
        \begin{aligned}
            \max_{{\bm{p}}\in\Pcal} \qquad & {\bm{p}}^\top{\bm{d}}\\ 
            \st \qquad & {\bm{d}} = \hat{\bm{\alpha}}+\hat B{\bm{p}},\\
            &Ad \le \frac{{\bm{c}}^t}{T-t+1}\\ 
            &{\bm{d}} \ge 0,
        \end{aligned}
    \end{equation*}
    there exists constant $C_2$ such that 
    \begin{equation*}
        \dist{\hat {\bm{d}},D^{\pi,t}} \le C_2\kappa^{-1}\norm{B-\hat B}_2
    \end{equation*}
    holds for optimal solution $\hat{\bm{d}}$ of the second constrained programming problem and all $\hat B$ such that $\norm{B-\hat B}_2<\delta$, where $\delta>0$ depends on $\lambda_{\min}(B+B^\top),$ the minimal eigenvalue of $B+B^\top$.
\end{lemma}
Note that, by Lemma \ref{lem:Lip1}, we only need to show that, there exists $\kappa>0$ such that $r({\bm{d}}) \le r({\bm{d}}^{\pi,t}) + \kappa(\dist{{\bm{d}},D^{\pi,t}})^2$ for any $d$ in the feasible set. Fortunately, the following lemma gives the existence of such constant.
\begin{lemma}
    \label{lem:second_order1}
    For any $d$ in the feasible set, we have $r({\bm{d}})-r({\bm{d}}^{\pi,t}) \ge \frac{1}{4}\lambda_{\min}(B+B^\top)\dist{{\bm{d}},{\bm{d}}^{\pi,t}}$.
\end{lemma}
By combining Lemma \ref{lem:Lip1} and \ref{lem:second_order1}, we now have 
\begin{align}\label{eq:bounded_by_params}
    \norm{\tilde {\bm{d}}^k-{\bm{d}}^{\pi,kn+1}}_2\le 4C_2\lambda_{\min}^{-1}(B+B^\top)\norm{\hat B^{kn+1}-B}_2,
\end{align}
By setting $\lambda = \delta$ in \eqref{eq:estimation_error_parameter1}, we get 
\begin{align*}
    \P(\norm{\hat b^{kn+1}-B}_2>\delta)\le C_1(kn)^{n^2+n-1}\exp\prn{-\frac{\sigma_0^2\sigma_1\sqrt{kn}}{8n}(\delta\wedge\delta^2)}.
\end{align*}
With a similar argument as above and in the proof of Theorem 6 in \citealt{keskin2014dynamic}, we can get:
\begin{align*}
    \ex{}{\norm{\tilde {\bm{d}}^k-{\bm{d}}^{\pi,kn+1}}_2^2}\le \frac{C_3n^5\lambda_{\min}^{-1}(B+B^\top)\max\{\delta^{-1},\delta^{-2}\}}{\sigma_1\sqrt{kn}},
\end{align*}
where $C_3$ is some constant determined by $C_1,C_2$. It follows that 
\begin{align}
    \label{eq:error_delta_3}
    \ex{}{\norm{\Delta_I^t}_2^2} &= \ex{}{\norm{\tilde {\bm{d}}^{k} -{\bm{d}}^{\pi,kn+1}}_2^2}\nonumber\\
    &\le  \frac{C_4n^5\lambda_{\min}^{-2}(B+B^\top)\max\{\delta^{-1},\delta^{-2}\}}{\sigma_1\sqrt{kn}}
\end{align}
for some constant $C_4$.
Plugging \eqref{eq:error_delta_1}, \eqref{eq:error_delta_2} and \eqref{eq:error_delta_3} into \eqref{eq:error_decomp} yields:
\begin{align}
    \label{eq:bound_delta}
    \ex{}{\norm{\Delta^t}_2^2}&\le 3\prn{\norm{\Delta_I^t}_2^2+\norm{\Delta_{II}^t}_2^2+\norm{\Delta_{III}^t}_2^2}\nonumber\\ 
    &\le 3\prn{\frac{4n^2d_{\max}^2}{(t-1)}+\frac{\sigma_0^2\norm{B}_2^2}{\sqrt{t}}+\frac{C_3n^5\lambda_{\min}^{-1}(B+B^\top)\max\{\delta^{-1},\delta^{-2}\}}{\sigma_1\sqrt{kn}}}\nonumber\\
    &\le \frac{C_5}{\sqrt{k}},
\end{align}
where $C_5 = 12\max\{4n^2d_{\max}^2,\sigma_0^2\norm{B}_2^2,C_3\sigma_1^{-1}n^{9/2}\lambda_{\min}^{-1}(B+B^\top)\max\{\delta^{-2},\delta^{-1}\}$.

Now we proceed the argument of \eqref{eq:regret_decomp2}. With a similar argument of \eqref{eq:single_step1}, we can get
\begin{align}
    &\label{single_Step3}
    \Rcal^T(\mix^{k},\F^T)-\Rcal^T(\mix^{k+1},\F^T)\nonumber\\
    &= (T-kn)r({\bm{d}}^{\pi,k})-\sum_{t=kn+1}^{k+1}r({\bm{d}}^{t},{\bm{\epsilon}}^t)-(T-(k+1)n)r({\bm{d}}^{\pi,(k+1)n}).
\end{align}
Now we proceed with the following three cases as similar as in Appendix \ref{appendix:resolve}: Case (I) $\min_i \tilde d_i^{t}\ge\zeta\brk{(T-t+1)^{-1/4}+t^{-1/4}},$ $\forall kn+1\le t\le(k+1)n$; Case (II):  $\max_i \tilde d_i^{t}<\zeta\brk{(T-t+1)^{-1/4}+t^{-1/4}},\forall kn+1\le t\le(k+1)n,\quad\forall kn+1\le t\le(k+1)n$; Case (III): others.

\paragraph*{Case (I): $\min_i \tilde p_i^{t}\ge\zeta\brk{(T-t+1)^{-1/4}+t^{-1/4}},$ $\forall kn+1\le t\le(k+1)n$.} We follow the streamline in Appendix \ref{appendix:resolve}, except we have different rounding threshold as well as the estimation error. We define $\Ecal_i^t = \{d_i^{\pi,t}\ge \frac{\epsilon_i^t+\Delta_i^t}{T-t}\}$. and $\Ecal^k=\cap_{i=1}^n\cap_{t=kn+1}^{(k+1)n}\Ecal_i^t.$

Recall the decomposition:
\begin{align}
    \label{eq:decomp20}
    &\ex{}{\Rcal^{T}(\mix^{kn},\F^T)-\Rcal^T(\mix^{k+1},\F^T)} 
    \nonumber\\&= \P(\Ecal^t)\ex{}{\Rcal^T(\mix^{t},\F^T)-\Rcal^T(\mix^{t+1},\F^T)\big\vert\Ecal^t}\nonumber\\&\qquad+ \P((\Ecal^t)^c)\ex{}{\Rcal^T(\mix^{t},\F^T)-\Rcal^T(\mix^{t+1},\F^T)\big\vert(\Ecal^t)^c}
\end{align}
With a similar argument of deriving \eqref{eq:decomp11},we have
\begin{align}
    \label{eq:decomp21}
    &\ex{}{\Rcal^{T}(\mix^{kn},\F^T)\Rcal^T(\mix^{(k+1)n},\F^T)\vert \Ecal^k}\nonumber\\ 
    &= \ex{}{(T-kn)r({\bm{d}}^{\pi,k})-\sum_{t=kn+1}^{(k+1)n}r({\bm{d}}^{t},{\bm{{\bm{\epsilon}}}}^t)-(T-(k+1)n)r({\bm{d}}^{\pi,k+1})\vert \Ecal^k}\nonumber\\
    &\le\ex{}{-\sum_{t=kn+1}^{(k+1)n}({\bm{\epsilon}}^t)^\top {\bm{p}}^t - \sum_{t=kn+1}^{(k+1)n}\brk{(\Delta^t)^\top B^{-1}\Delta^t+(2{\bm{d}}^{\pi,k}-{\bm{\alpha}})^\top B^{-1}\Delta^t}\vert\Ecal^k}\nonumber\\ 
    &\qquad+(T-(k+1)n)\ex{}{\frac{(\sum_{t=kn+1}^{(k+1)n}{\bm{\epsilon}}^t+\Delta^t)^\top}{T-(k+1)n}B^{-1}\frac{\sum_{t=kn+1}^{(k+1)n}({\bm{\epsilon}}^t+\Delta^t)}{T-(k+1)n}+(2{\bm{d}}^{\pi,k}-{\bm{\alpha}})^\top B^{-1}\frac{\sum_{t=kn+1}^{(k+1)n}\Delta^t+{\bm{\epsilon}}^t}{T-(k+1)n}\vert\Ecal^k}\nonumber\\
    & =\ex{}{({\bm{d}}^{\pi,k})^\top B^{-1}\sum_{t=kn+1}^{(k+1)n}{\bm{\epsilon}}^t-\sum_{t=kn+1}^{(k+1)n}({\bm{\epsilon}}^t)^\top B^{-1}\Delta^t\vert\Ecal^k} + \frac{1}{T-(k+1)n}\ex{}{(\sum_{t=kn+1}^{(k+1)n}{\bm{\epsilon}}^t+\Delta^t)^\top B^{-1}(\sum_{t=kn+1}^{(k+1)n}{\bm{\epsilon}}^t+\Delta^t)\vert\Ecal^k}.
\end{align}
Note that ${\bm{\epsilon}}^t$ is independent of $\Delta^t.$ We have $\ex{}{{\bm{\epsilon}}^t B^{-1}\Delta^t} = 0.$ Note that $\Delta^t\le {\bm{d}}_{\max}.$ Following the proof of \eqref{eq:decomp12} and \eqref{eq:decomp13} leads to 
\begin{align}
\label{eq:decomp22}
    &\P(\Ecal^k)\ex{}{({\bm{d}}^{\pi,k})^\top B^{-1}\sum_{t=kn+1}^{(k+1)n}{\bm{\epsilon}}^t-\sum_{t=kn+1}^{(k+1)n}({\bm{\epsilon}}^t)^\top B^{-1}\Delta^t\vert\Ecal^k}\nonumber\\ 
    &\le 2\norm{B^{-1}}_2d_{\max}\P((\Ecal^k)^c)^{1/2}\ex{}{\sum_{t=kn+1}^{(k+1)n}\norm{{\bm{\epsilon}}^t}_2^2}^{1/2}\nonumber\\
    &\le 2\sigma\norm{B^{-1}}_2d_{\max}\P((\Ecal^k)^c)^{1/2}
\end{align}
and 
\begin{align}
    \label{eq:decomp23}
    &\P(\Ecal^k)\ex{}{(\sum_{t=kn+1}^{(k+1)n}{\bm{\epsilon}}^t+\Delta^t)^\top B^{-1}(\sum_{t=kn+1}^{(k+1)n}{\bm{\epsilon}}^t+\Delta^t)\vert\Ecal^k}\nonumber\\
    &\le \frac{2n\norm{B^{-1}}_2}{T-(k+1)n}\sum_{t=kn+1}^{(k+1)n}\ex{}{\norm{{\bm{\epsilon}}^t}_2^2+\norm{\Delta^t}_2^2}\nonumber\\
    &\overset{(a)}{\le} \frac{2n\norm{B^{-1}}_2}{T-(k+1)n}\sum_{t=kn+1}^{(k+1)n}\prn{\sigma^2+C_5k^{-1/2}}\nonumber\\
    &= \frac{2n^2\norm{B^{-1}}_2}{T-(k+1)n}\prn{\sigma^2+C_5k^{-1/2}},
\end{align}
where we apply \eqref{eq:bound_delta} in (a).

Now following the deduction in \eqref{eq:decomp14}, we get
\begin{align}
\label{eq:decomp24}
    &\P((\Ecal^t)^c)\ex{}{\Rcal^T(\mix^{t},\F^T)-\Rcal^T(\mix^{t+1},\F^T)\big\vert(\Ecal^t)^c}\nonumber\\
    &\le (T-kn)r_{\max}\P((\Ecal^k)^c)+\sqrt{n}\sigma\P((\Ecal^k)^c)^{1/2}.
\end{align}
Now combining \eqref{eq:decomp20}, \eqref{eq:decomp21}, \eqref{eq:decomp22}, \eqref{eq:decomp23} and \eqref{eq:decomp24}, we get 
\begin{align}\label{eq:single_step222}
    &\ex{}{\Rcal^T(\mix^{t},\F^T)-\Rcal^T(\mix^{t+1},\F^T)}\nonumber\\
   & 2\sigma\norm{B^{-1}}_2d_{\max}\P((\Ecal^k)^c)^{1/2} + \frac{2n^2\norm{B^{-1}}_2}{T-(k+1)n}(\sigma^2+C_5k^{-1/2})+(T-kn)r_{\max}\P((\Ecal^k)^c)+\sqrt{n}\sigma\P((\Ecal^k)^c)^{1/2}.
\end{align}
Now we give bound to $\P((\Ecal^t)^c).$ We consider two events: $\Ical_{i1}^t = \{\zeta((T-t+1)^{-1/4}+t^{-1/4})/3 \ge |\Delta_i^t|\};$ $\Ical_{i2}^t = \{\zeta((T-t+1)^{-1/4}+t^{-1/4})/3\ge \frac{\epsilon_i^t}{T-(k+1)n}\}$. Then we have 
\begin{align*}
    \P((\Ecal^k)^c)&\le 1-\P(\cap_{i=1}^n\cap_{t=kn+1}^{(k+1)n}(\Ical_{i1}^t\cap\Ical_{i2}^t)) \\
    &\le \sum_{i=1}^n\sum_{t=kn+1}^{(k+1)n}(\P(\Ical_{i1}^t)+\P(\Ical_{i2}^t)).
\end{align*}
With \eqref{eq:bounded_by_params}, \eqref{eq:error_delta_1}, \eqref{eq:error_delta_2} and \eqref{eq:estimation_error_parameter1}, we have 
\begin{align}
    \label{eq:Ical1}
    &\P((\Ical_{i1}^t)^c)\nonumber\\
    &\le \P\prn{4C_2\lambda_{\min}^{-1}(B+B^\top)\norm{\hat B^{kn+1}-B}_2+\frac{2n}{t-1}{\bm{d}}_{\max}+\sigma_0\norm{B}_2t^{-1/4} > \zeta((T-t+1)^{-1/4}+t^{-1/4})/3}\nonumber\\ 
    &\le \P\prn{C_2\lambda_{\min}^{-1}(B+B^\top)\norm{\hat B^{kn+1}-B}_2 > \zeta((T-t+1)^{-1/4}+t^{-1/4})/24}\nonumber\\ 
    &\le C_1(kn)^{n^2+n-1}\exp\prn{-\frac{\sigma_0^2\sigma_1\sqrt{kn}}{8n}(\lambda\wedge\lambda^2)}\nonumber\\ 
    &\le C_5^2/(n^2T^2)
\end{align}
for some constant $C_5$, where $\lambda = C_2^{-1}\lambda_{\min}(B+B^\top)\zeta((T-t+1)^{-1/4}+t^{-1/4})/24$ and we have 
\begin{align*}
    \lambda^2\wedge\lambda\ge \frac{8(n^{5/2}+n)\log^3 n\log^2 T}{k^{-1/2}\sigma_0^2\sigma}
\end{align*}
by definition of $\zeta.$ For the event $\Ical_{i2}^c,$ it follows directly from the derivation of \eqref{eq:prob_complementary} that 
\begin{align}
    \label{eq:Ical2}
    &\P((\Ical_{i2}^t)^c\nonumber\\ 
    &\le \P\prn{\zeta((T-t+1)^{-1/4}+t^{-1/4})/3< \frac{\epsilon_i^t}{T-(k+1)n}}\nonumber\\
    &\le \exp\prn{-\frac{\zeta^2(T-(k+1)n)^2}{36\sigma^2}}\nonumber\\ 
    &\le C_5^2/(n^2T^2).
\end{align}
Now combining \eqref{eq:Ical1} and \eqref{eq:Ical2}, we get
\begin{align*}
     \P((\Ecal^k)^c)&\le \sum_{i=1}^n\sum_{t=kn+1}^{(k+1)n}(\P(\Ical_{i1}^t)+\P(\Ical_{i2}^t))\\
     &\le 2C_5^2/T^2.
\end{align*}
Plugging the above inequality into \eqref{eq:single_step222} yields
\begin{align}
    \label{eq:case_i_summary}
        &\ex{}{\Rcal^T(\mix^{k},\F^T)-\Rcal^T(\mix^{k+1},\F^T)}\nonumber\\
   & \le 2C_5\sigma\norm{B^{-1}}_2d_{\max}T^{-1} + \frac{2n^2\norm{B^{-1}}_2}{T-(k+1)n}(\sigma^2+C_5k^{-1/2})+C_5^2(T-kn)r_{\max}/T^2+C_5\sqrt{n}\sigma T^{-1}.
\end{align}

\paragraph*{Case (II):  $\max_i \tilde d_i^{t}<\zeta\brk{(T-t+1)^{-1/4}+t^{-1/4}},\forall kn+1\le t\le(k+1)n,\quad\forall kn+1\le t\le(k+1)n$.} This case is much simpler, just follow \eqref{eq:single_step22}:
\begin{align}
    \label{eq:case_ii_summary}
        &\ex{}{\Rcal^T(\mix^{k},\F^T)-\Rcal^{T}(\mix^{k+1},\F^T)}\nonumber\\ 
        &\le \ex{}{(T-kn)r({\bm{d}}^{\pi,k})-(T-(k+1)n)r\prn{\frac{T-kn}{T-(k+1)n}{\bm{d}}^{\pi,k}} }\nonumber\\ 
        &= \ex{}{(T-kn)({\bm{d}}^{\pi,k})^\top B^{-1}({\bm{d}}^{\pi,t}-{\bm{\alpha}})- (T-(k+1)n)\frac{T-kn}{T-(k+1)n}({\bm{d}}^{\pi,t})^\top B^{-1}(\frac{T-kn}{T-(k+1)n}{\bm{d}}^{\pi,t}-{\bm{\alpha}} ) }\nonumber\\ 
        &= \ex{}{-\frac{T-kn}{T-(k+1)n}({\bm{d}}^{\pi,t})^\top B^{-1}{\bm{d}}^{\pi,t}}\nonumber\\
        &=-\frac{T-kn}{T-(k+1)n}\ex{}{({\bm{d}}^{\pi,t}+\Delta^t)^\top B^{-1}({\bm{d}}^t+\Delta^)}\nonumber\\
        &\le -\frac{2(T-kn)}{T-(k+1)n}\ex{}{({\bm{d}}^{kn+1})^\top B^{-1}{\bm{d}}^{kn+1}}-\frac{2(T-kn)}{T-(k+1)n}\ex{}{(\Delta^{kn+1})^\top B^{-1}\Delta^{kn+1}}\nonumber\\
        &\le 8\zeta^2\brk{(T-kn+1)^{-1/2}+(kn)^{-1/2}}+\frac{8C_5\norm{B^{-1}}_2}{\sqrt{k}}.
\end{align}
where we again use \eqref{eq:bound_delta} in the last line.

\paragraph*{Case (III): others.} Following the same streamline of the argument in Appendix \ref{appendix:resolve} and utilizing \eqref{eq:case_i_summary}, \eqref{eq:case_ii_summary}, we can get 
\begin{align}
    \label{eq:case_iii_summary}
    &\ex{}{\Rcal^T(\mix^{k},\F^T)-\Rcal^{T}(\mix^{k+1},\F^T)}\nonumber\\ 
    &\le 8\zeta^2\brk{(T-kn+1)^{-1/2}+(kn)^{-1/2}}+\frac{8C_5\norm{B^{-1}}_2}{\sqrt{k}}\nonumber\\
    &\qquad+2C_5\sigma\norm{B^{-1}}_2d_{\max}T^{-1} + \frac{2n^2\norm{B^{-1}}_2}{T-(k+1)n}(\sigma^2+C_5k^{-1/2})+C_5^2(T-kn)r_{\max}/T^2+C_5\sqrt{n}\sigma T^{-1}.
\end{align}

\paragraph*{Wrap-up.} Combinine \eqref{eq:case_i_summary}, \eqref{eq:case_ii_summary} and \eqref{eq:case_iii_summary} leads to 
\begin{align*}
    &\ex{}{\Rcal^T(\mix^{k},\F^T)-\Rcal^{T}(\mix^{k+1},\F^T)}\nonumber\\ 
    &\le 8\zeta^2\brk{(T-kn+1)^{-1/2}+(kn)^{-1/2}}+\frac{8C_5\norm{B^{-1}}_2}{\sqrt{k}}\nonumber\\
    &\qquad+2C_5\sigma\norm{B^{-1}}_2d_{\max}T^{-1} + \frac{2n^2\norm{B^{-1}}_2}{T-(k+1)n}(\sigma^2+C_5k^{-1/2})+C_5^2(T-kn)r_{\max}/T^2+C_5\sqrt{n}\sigma T^{-1}.
\end{align*}
Plugging this into \eqref{eq:decomp20} yields
\begin{align*}
    \regret[T]{\pi} &=\sum_{k=1}^{T'}\ex{}{\Rcal^T(\mix^{k},\F^T)-\Rcal^T(\mix^{k+1},\F^T)}\\ 
    &=O(\brk{(\zeta^2+C_5\norm{B^{-1}}_2}\sqrt{T}).
\end{align*}

\subsection{Proof of Proposition \ref{prop:impossible}}
We follow Sec 15.2 \cite{lattimore2020bandit} and \cite{cheung2024leveraging}. Consider two demand functions ${\bm{d}}_{\theta}({\bm{p}}) = 2\Delta-\Delta{\bm{p}}+{\bm{\epsilon}}$ and ${\bm{d}}_{\theta'}({\bm{p}}) = 3\Delta-2\Delta{\bm{p}}+{\bm{\epsilon}}$ with $\Delta = T^{-\beta},$ where ${\bm{\epsilon}}\sim\N(0,1).$ Assume that both $\theta,\theta'$ have same offline data distribution $\P^\off$. Denote $P_\pi^{\on},Q_\pi^{\on}$ as the online distribution of the demands under policy $\pi$. Additionally, define $P,Q$ as the joint distribution of $(P^\off,P_\pi^\on),(P^\off,Q_\pi^\on),$ respectively. Then $\regret{\pi,\theta} = \ex{P}{\sum_{t=1}^T\Delta^{-1}(\Delta{\bm{p}}^t-\Delta)^2}$ and $\regret{\pi,\theta'} = \ex{Q}{\sum_{t=1}^T2\Delta^{-1}(\Delta{\bm{p}}^t-\frac{3\Delta}{2})^2}$. For set $S\subset R$, denote $N_T(S)$ as the total number of $p_t$ such that $p_t\in S$. Consider the event $\Ical=\{N_T([\frac{5\Delta}{4},\infty))>\frac{T}{2}\}$. Then by Bretagnolle-Huber inequality \citealp[Thm 14.2]{lattimore2020bandit} we have
\begin{align*}
    &\regret{\pi,\theta}+\regret{\pi,\theta'}\\
    &\ge  \frac{\Delta}{16}*\frac{T}{2}*P(\Ical)+\frac{2\Delta}{16}*\frac{T}{2}*Q(\Ical^C)\\
    &\ge \frac{T^{1-\gamma}}{32}(P(\Ical)+Q(\Ical^C))\\
    &\ge \frac{T^{1-\gamma}}{32}\exp(-\KL{P}{Q}).
\end{align*}
For Gaussian $P_\theta({\bm{p}})=\N(2\Delta-\Delta{\bm{p}},1)$ and $Q_{\theta'}({\bm{p}})=\N(3\Delta-2\Delta{\bm{p}},1),$ we have $\KL{P_\theta}{P_{\theta'}}=\frac{(\Delta-p)^2}{2}.$ In order to to further give bound to the term $\exp(-\KL{P}{Q})$, we introduce the following lemma in \cite{cheung2024leveraging}.
\begin{lemma}
    \label{lem:offline}
    Consider two instances with onlien distribution $P$, $Q$ and shared offline dataset with samples $\{(p^{-t},{\bm{d}}^{-t})\}_{t=1}^N$, then for any admissible policy $\pi$, it holds that 
    \begin{align*}
        \exp(-\KL{P}{Q})&=\exp(-\ex{P}{\sum_{t=1}^T\KL{P_\theta(p^t)}{Q_{\theta'}(p^t)} }.
    \end{align*}
\end{lemma}
By Lemma \ref{lem:offline}, we have
\begin{align*}
    \exp(-\KL{P}{Q})&=\exp(-\ex{P}{\sum_{t=1}^T\KL{P_\theta(p^t)}{Q_{\theta'}(p^t)}}\\
    &=\exp(-\ex{P}{\sum_{t=1}^T\frac{(\Delta-\Delta{\bm{p}}^t)^2}{2}}\\ 
    & = \exp(-\Delta\regret{\pi,\theta})\\ 
    &\ge \exp(-\Delta CT^\gamma)\\ 
    & = \exp(-C).
\end{align*}
As a result, we have 
\begin{align*}
    \regret{\pi,\theta}+\regret{\pi,\theta'}\ge \frac{T^{1-\gamma}}{32}\exp(-C)=\Omega(T^{1-\gamma}).
\end{align*}
Since $\beta\in[0,\frac{1}{2}),$ we know that $\regret{\pi,\theta'}\ge \Omega(T^{1-\gamma}).$

\section{Proof of Theorem \ref{thm:incumbent}}\label{appendix:incumbent}
In this section, we follow a similar streamline as in Appendix \ref{appendix:learn}. We only need to consider the case $({\bm{\epsilon}}^0)^2T\le \rho\sqrt{T}$ and the informed prices are utilized. To begin with, recall our re-solve constrained programming problem:
\begin{equation}
    \label{prob:resolve_restate3}
    \begin{aligned}
        \max_{p\in\Pcal}\qquad & r ={\bm{p}}^\top{\bm{d}}\\ 
    \st\qquad& {\bm{d}} = {\bm{\alpha}} + B{\bm{p}},\\ 
    &Ad\le \frac{{\bm{c}}^t}{T-t+1},
    \end{aligned}
\end{equation}
where ${\bm{c}}^t$ is the inventory level at the beginning of time $t$. We use the single-step difference decomposition:
\begin{align}
    \label{eq:regret_decomp3}
    \regret[T]{\pi} &= \ex{}{\sum_{t=1}^{T}\Rcal^T(\mix^{t},\F^T)-\Rcal^T(\mix^{t+1},\F^T)}\nonumber\\
    &=\sum_{t=1}^{T}\ex{}{\Rcal^T(\mix^{t},\F^T)-\Rcal^T(\mix^{t+1},\F^T)}.
\end{align}
We proceed by giving bound to the term $\ex{}{\Rcal^T(\mix^{t},\F^T)-\Rcal^T(\mix^{t+1},\F^T)}$ for $1\le t\le T.$ We define $\Delta^t =  {\bm{d}}^t - {\bm{d}}^{\pi,t}$ similarly.  Let $l=\mod(t,n)$. Then we have 
\begin{align}\label{eq:error_decomp2}
   \Delta^t = \underset{:=\Delta_I^t}{\underbrace{(\tilde {\bm{d}}^t-{\bm{d}}^{t})}} + \underset{:=\Delta_{II}^t}{\underbrace{\sigma_0t^{-1/2}Be_{l}}}.
\end{align}
It follows directly that 
\begin{align}
    \label{eq:error_delta_21}
\norm{\Delta_{II}^t}_2\le\sigma_0t^{-1/2}\norm{B}_2.
\end{align}
In order to give bound to the term $\Delta_I^t,$ we define the quanity $J^t$ as 
\begin{align}
    \label{eq:J_t2}
    J^t := n^{-1}\sum_{s=1}^t\norm{p^s-p^0}_2^2\ge ^{-1}\sum_{s=1}^t\norm{p^\star-p^0}_2^2=tn^{-1}\delta_0^2.
\end{align}
Now we introduce the following tail bound on the parameter estimation analogous to Lemma \ref{lem:tail}.
\begin{lemma}[\citealt{keskin2014dynamic}]
\label{lem:tail2}
    Under our policy and estimation \eqref{eq:regression_informed}, if ${\bm{\epsilon}}^0=0$ (i.e. no misspecification), then there exist finite positive constants $C_6,\sigma_2$ such that
    \begin{align*}
        \P(\norm{\hat B^t-B}_2>\lambda,J^t\ge\lambda')\le C_6t^{n^2-1}\exp(-\sigma_2(\lambda\wedge\lambda^2)J^t).
    \end{align*}
\end{lemma}
Moreover, we have the following lower bound on the minimal eigenvalue of the matrix $\hat P^t:=\sum_{s=1}^{t-1}(p^s-p^0)(p^s-p^0)^\top$ in \eqref{eq:regression_informed}:
\begin{lemma}[\citealt{keskin2014dynamic}]
\label{lem:minimal_eigenvalue}
    Under our policy, we have $\lambda_{\min}(\hat P^t)\ge J^t.$
\end{lemma}
Note that $\hat P^t$ is positive definite by the deduction above. Combining \eqref{eq:regression_informed}, Lemma \ref{lem:tail2} and \ref{lem:minimal_eigenvalue} leads to 
\begin{align}
\label{eq:tail}
    \P(\norm{\hat B^t-B}_2>\lambda + {\bm{\epsilon}}^0\le C_7t^{n^2-1}\exp(-\sigma_2(\lambda\wedge\lambda^2)J^t)
\end{align}
for some constant $C_7$ under offline data misspecification ${\bm{\epsilon}}^0$. Now following a similar argument as in the proof of \eqref{eq:error_delta_3} (except that we have the $J^t$ linear in $T$, rather than square root order of $T$ in the previous section, hence we have better bound here), we can get 
\begin{align}
    \label{eq:error_delta_22}
    \ex{}{\norm{\Delta_I^t}_2^2}\le \frac{C_8n^5\lambda_{\min}^{-1}(B+B^\top)\max\{\delta^{-1},\delta^{-2}\}}{\sigma_2t}+({\bm{\epsilon}}^0)^2.
\end{align}
Plugging \eqref{eq:error_delta_21} and \eqref{eq:error_delta_22} into \eqref{eq:error_decomp2}
\begin{align}
    \label{eq:bound_delta2}
    \ex{}{\norm{\Delta^t}_2^2}&\le 2\prn{\norm{\Delta_I^t}_2^2+\norm{\Delta_{II}}_2^2}\nonumber\\
    &\le 2\prn{\frac{C_8n^5\lambda_{\min}^{-1}(B+B^\top)\max\{\delta^{-1},\delta^{-2}\}}{\sigma_2t}+({\bm{\epsilon}}^0)^2+\sigma_0^2t^{-1}\norm{B}_2^2}\nonumber\\ 
    &\le \frac{C_9}{t},
\end{align}
where $C_9 = 8\max\{C_8n^5\lambda_{\min}^{-1}(B+B^\top)\max\{\delta^{-1},\delta^{-2}\}/\sigma_2,\sigma_0^2\norm{B}_2^2\}$.

Now we proceed with almost the same argument as in Appendix \ref{appendix:learn}. We consider three cases: Case (I) $\min_i\tilde d_i^t\ge \zeta\brk{(T-t+1)^{-1/2}+t^{-1/2}}$; Case (II) $\max_i\tilde d_i^t\le \zeta\brk{(T-t+1)^{-1/2}+t^{-1/2}}$; Case (III) others.

\paragraph*{Case (I).} We follow a almost same argument as the derivation of \eqref{eq:case_i_summary} except by replacing the bound of $\ex{}{\norm{\Delta}_2^2}$ by \eqref{eq:bound_delta2} and the tail bound of \eqref{eq:Ical1}, \eqref{eq:Ical2} with corresponding rounding threshold $\zeta\brk{(T-t+1)^{-1/2}+t^{-1/2}}$ and $J^t$ in \eqref{eq:J_t2}. We can get
\begin{align}
    \label{eq:case_i_summary2}
    &\ex{}{\Rcal^T(\mix^t,\F^T)-\Rcal^T(\mix^{t+1},\F^T)}\nonumber\\
    &4C_9\sigma_0\norm{B^{-1}}{\bm{d}}_{\max}T^{-1}+\frac{4n^2\norm{B^{-1}}_2}{T-t}(\sigma^2+C_9T^{-1})+2C_9^2(T-t+1)r_{\max}/T^2+C_0\sqrt{n}\sigma T^{-1}+2({\bm{\epsilon}}^0)^2.
\end{align}

\paragraph*{Case (II).} With a similar argument of deriving \eqref{eq:case_ii_summary} by replacing the rounding threshold with $\zeta\brk{(T-t+1)^{-1/2}+t^{-1/2}}$, we can get 
\begin{align}
    \label{eq:case_ii_summary2}
    &\ex{}{\Rcal^T(\mix^t,\F^T)-\Rcal^T(\mix^{t+1},\F^T)}\nonumber\\
    &\le 8\zeta^2\brk{(T-t+1)^{-1}+t^{-1}}+\frac{8C_9\norm{B^{-1}}_2}{t}.
\end{align}

\paragraph*{Case (III).} With a similar argument of deriving \eqref{eq:case_iii_summary2}, we arrive at 
\begin{align}
    \label{eq:case_iii_summary2}
    &\ex{}{\Rcal^T(\mix^t,\F^T)-\Rcal^T(\mix^{t+1},\F^T)}\nonumber\\
    &\le 4C_9\sigma_0\norm{B^{-1}}{\bm{d}}_{\max}T^{-1}+\frac{4n^2\norm{B^{-1}}_2}{T-t}(\sigma^2+C_9T^{-1})+2C_9^2(T-t+1)r_{\max}/T^2+C_0\sqrt{n}\sigma T^{-1}+2({\bm{\epsilon}}^0)^2\nonumber\\ 
    &\qquad+8\zeta^2\brk{(T-t+1)^{-1}+t^{-1}}+\frac{8C_9\norm{B^{-1}}_2}{t}.
\end{align}

\paragraph*{Wrap-up.} Now combining \eqref{eq:case_i_summary2}, \eqref{eq:case_ii_summary2} and \eqref{eq:case_iii_summary2} yields 
\begin{align*}
    &\ex{}{\Rcal^T(\mix^t,\F^T)-\Rcal^T(\mix^{t+1},\F^T)}\nonumber\\
        &\le 4C_9\sigma_0\norm{B^{-1}}{\bm{d}}_{\max}T^{-1}+\frac{4n^2\norm{B^{-1}}_2}{T-t}(\sigma^2+C_9T^{-1})+2C_9^2(T-t+1)r_{\max}/T^2+C_0\sqrt{n}\sigma T^{-1}+2({\bm{\epsilon}}^0)^2\nonumber\\ 
    &\qquad+8\zeta^2\brk{(T-t+1)^{-1}+t^{-1}}+\frac{8C_9\norm{B^{-1}}_2}{t}.
\end{align*}
Now plugging the above inequality into \eqref{eq:regret_decomp3} leads to
\begin{align*}
        \regret[T]{\pi} 
    &=\sum_{t=1}^{T}\ex{}{\Rcal^T(\mix^{t},\F^T)-\Rcal^T(\mix^{t+1},\F^T)}\nonumber\\ 
    &= O\prn{C_{10}\log T+({\bm{\epsilon}}^0)^2T},
\end{align*}
where $C_{10} = 2C_0\sigma
_0\norm{B^{-1}}_2d_{\max}+2n^2\norm{B^{-1}}_2\sigma
^2+C_0\sqrt{n}\sigma+16\zeta^2+8C_9\norm{B^{-1}}_2.$

\section{Proof of Proposition \ref{prop:impossible_incumbent}}\label{appendix:impossible_incumbent}
This is a direct result following Proposition \ref{prop:impossible}. We just need to take ${\bm{\epsilon}}$ of order $T^{-(1-\gamma)/2}$.



\section{Proof of lemmas}

\subsection{Proof of Lemma \ref{lem:second_order1}}
We rewrite the target as $r({\bm{d}}) = {\bm{d}}^\top B^{-1}({\bm{d}}-{\bm{\alpha}},$ which is a quadratic function since $B$ is positive definite. For optimal solution ${\bm{d}}^{\pi,t},$ ${\bm{d}}-{\bm{d}}^{\pi,t}$ is an descending direction for feasible $d$ and it follows from \cite{boyd2004convex} that $({\bm{d}}^\top-{\bm{d}}^{\pi,t})\nabla_{{\bm{d}}^{\pi,t}}r({\bm{d}}^{\pi,t})\le 0.$ Moreover, with direct calculation, we can get
\begin{align*}
    r({\bm{d}}) = r({\bm{d}}^{\pi,t})+ ({\bm{d}}-{\bm{d}}^{\pi,t})^\top\nabla_{{\bm{d}}^{\pi,t}}r({\bm{d}}^{\pi,t}) + \frac{1}{2}({\bm{d}}-{\bm{d}}^{\pi,t})^\top B^{-1}({\bm{d}}-{\bm{d}}^{\pi,t}).
\end{align*}
As a result, we have
\begin{align*}
    r({\bm{d}}) &\le r({\bm{d}}^{\pi,t}) + \frac{1}{2}({\bm{d}}-{\bm{d}}^{\pi,t})^\top B^{-1}({\bm{d}}-{\bm{d}}^{\pi,t})\\
   & \le r({\bm{d}}^{\pi,t})+\frac{1}{4}\lambda_{\min}(B^{-1}+B^{-\top})\norm{{\bm{d}}-{\bm{d}}^{\pi,t}}_2^2\\
   &\le r({\bm{d}}^{\pi,t})+\frac{1}{4}\lambda_{\min}(B^{-1}+B^{-\top}).
\end{align*}

\end{document}